\magnification1200

\centerline{\bf The dense amalgam of metric compacta} 
\centerline{\bf and topological characterization of
boundaries of free products of groups}

\bigskip
\centerline{Jacek \'Swi\c atkowski
\footnote{*}{The author was partially supported by the Polish National
Science Centre (NCN), Grant 2012/06/ST1/00259.}}

\bigskip
\centerline{Instytut Matematyczny, Uniwersytet Wroc\l awski}
\centerline{pl. Grunwaldzki 2/4, 50-384 Wroc\l aw, Poland}
\centerline{E-mail: {\tt swiatkow@math.uni.wroc.pl}}

\bigskip
{\bf Abstract.} We introduce and study the operation, called 
{\it dense amalgam}, which to any tuple $X_1,\dots,X_k$ of
non-empty compact metric spaces associates some disconnected
perfect compact metric space, denoted $\widetilde\sqcup(X_1,\dots,X_k)$,
in which there are many appropriately distributed copies of
the spaces $X_1,\dots,X_k$.
We then show that, in various settings, the ideal boundary of the free product 
of groups is homeomorphic to the dense amalgam of boundaries
of the factors. We give also related more general results for graphs of groups
with finite edge groups. We justify these results by referring to a convenient
characterization of dense amalgams, in terms of a list of properties,
which we also provide in the paper.
As another application, we show that the boundary of a Coxeter group
which has infinitely many ends, and which is not virtually free,
is the dense amalgam of the boundaries of its maximal 1-ended special
subgroups.

\bigskip\noindent
{\bf 0. Introduction.}

\smallskip
In Section 1 of the paper we describe an operation which to any finite tuple $X_1,\dots,X_k$ of nonempty metric compacta associates a  
metric compactum $$Y=\widetilde\sqcup(X_1,\dots,X_k)$$ 
which satisfies the following.
$Y$ can be equipped with a countable infinite family $\cal Y$ of subsets,
partitioned as ${\cal Y}={\cal Y}_1\sqcup\dots\sqcup{\cal Y}_k$,
such that:

\itemitem{
(a1)} the subsets in $\cal Y$ are pairwise disjoint and for each
$i\in\{ 1,\dots,k \}$ the family ${\cal Y}_i$ consists of embedded copies of the space $X_i$;

\itemitem{
(a2)} the family $\cal Y$  
is {\it null}, i.e. for any metric on $Y$ compatible with the topology the diameters of sets in $\cal Y$ converge to 0;

\itemitem{
(a3)} each $Z\in{\cal Y}$ is a boundary subset of $Y$ (i.e. its complement is dense);

\itemitem{
(a4)} for each $i$, the union of the family ${\cal Y}_i$ is dense in $Y$;

\itemitem{
(a5)} any two points of $Y$ which do not belong to the same subset
of $\cal Y$ can be separated from each other by an open and closed subset $H\i Y$ which is $\cal Y$-saturated (i.e. such that any element of $\cal Y$ is either contained in or disjoint with $H$).

\noindent
We call the operation $\widetilde\sqcup$ {\it the dense amalgam}, 
and its result $\widetilde\sqcup(X_1,\dots,X_k)$ 
{\it the dense amalgam of the spaces $X_1,\dots,X_k$}.  
Obviously, the dense amalgam of any family $X_1,\dots,X_k$
of spaces is a disconnected perfect compact metric space.
Moreover, if the spaces $X_1,\dots,X_k$ are connected, then the connected
components of their dense amalgam are the subsets from the family
$\cal Y$ and the singletons from the complement of the union $\cup{\cal Y}$.
In Section 3 we show that the operation 
of dense amalgam
satisfies
the following 
"algebraic"
properties.

\vfill\break

\medskip\noindent
{\bf 0.1 Proposition.} 

\item{(1)} {\it $\widetilde\sqcup(X_1,\dots,X_k)=
\widetilde\sqcup(X_1\sqcup\dots\sqcup X_k)$ 
(so in particular the operation is commutative).} 

\item{(2)} {\it $\widetilde\sqcup(X_1,\dots,X_k)=\widetilde\sqcup(X_1,\dots,X_{i-1},\widetilde\sqcup(X_i,\dots, X_k))$ for any $k\ge1$ and any $1\le i\le k$ (so the operation is associative and idempotent).}

\item{(3)} {\it $\widetilde\sqcup(X,X_1,\dots,X_k)=\widetilde\sqcup(X,X,X_1,\dots,X_k)$ for any $k\ge0$.}

\item{(4)} {\it For any totally disconnected nonempty compact metric space $Q$, and any $k\ge1$, we have}
$$
\widetilde\sqcup(X_1,\dots,X_k,Q)=\widetilde\sqcup(X_1,\dots,X_k).
$$

\item{(5)} {\it For any totally disconnected space $Q$ we have $\widetilde\sqcup(Q)=C$, where $C$ is the Cantor space.}

\medskip
In Section 2 we prove the following characterization result.

\medskip\noindent
{\bf 0.2 Theorem.}
{\it Given any nonempty compact metric spaces $X_1,\dots,X_k$,
each metric compactum $Y$ which can be equipped with a family 
${\cal Y}={\cal Y}_1\sqcup\dots\sqcup{\cal Y}_k$ 
of subsets satisfying conditions (a1)-(a5) above
is homeomorphic to the dense amalgam $\widetilde\sqcup(X_1,\dots,X_k)$.}

\medskip
The main motivation for the study of dense amalgams in this paper
comes from their role in understanding ideal boundaries of spaces and groups.
In Sections 4--6 we deal with various settings for
ideal boundaries,
showing among others that in these settings the boundary of the free product
of groups is homeomorphic to the dense amalgam of the boundaries of the factors. We also give similar more general 
results for fundamental groups of non-elementary 
graphs of groups with finite edge groups
(see Theorem 0.3 below). 
The term {\it non-elementary} for a graph of groups is explained in
Definition 4.1.5, but it contains among others the case when
the vertex groups are all infinite and the underlying graph is not
reduced to a single vertex. In particular, the results apply to
amalgamated free products and HNN extensions of infinite groups 
along finite subgroups.
In consistency with the properties from Proposition 0.1, we use the convention that $\widetilde\sqcup(\emptyset)$ is the Cantor space and $\widetilde\sqcup(\emptyset,X_1,\dots,X_k):=\widetilde\sqcup(X_1,\dots,X_k)$. The concepts of $E{\cal Z}$-boundaries appearing in the
statement below
are explained in Subsection 4.2 (see Definitions 4.2.1 and 4.2.2).
Our main result concerning boundaries is as follows.

\medskip\noindent
{\bf 0.3 Theorem.} {\it Let $G=\pi_1({\cal G})$, where $\cal G$ is 
a non-elementary graph of groups with finite edge groups. Let $v_1,\dots,v_k$ be the vertices of the underlying graph of $\cal G$.}
\item{(1)} {\it Suppose that all vertex groups $G_{v_i}$ of $\cal G$ admit  $E{\cal Z}$-boundaries (respectively, $E{\cal Z}$-boundaries in the strong sense of Carlsson-Pedersen), and let $\partial G_{v_i}$ be such boundaries. Then $G$ admits an $E{\cal Z}$-boundary (in the strong sense of Carlsson-Pedersen, respectively) homeomorphic to the dense amalgam $\widetilde\sqcup(\partial G_{v_1},\dots,\partial G_{v_k})$.}
\item{(2)} {\it Suppose that all vertex groups $G_{v_i}$ of $\cal G$ are hyperbolic and let $\partial G_{v_i}$ be their Gromov boundaries. Then the Gromov boundary $\partial G$ is homeomorphic to the dense amalgam $\widetilde\sqcup(\partial G_{v_1},\dots,\partial G_{v_k})$. }
\item{(3)} {\it Suppose that all vertex groups $G_{v_i}$ of $\cal G$ are CAT(0), and for each $v_i$ let $\Delta_{i}$ be a CAT(0) space on which $G_{v_i}$ acts geometrically. Then there is a CAT(0) space $\Delta$ on which $G$ acts geometrically, and such that the CAT(0) boundary $\partial\Delta$  is homeomorphic to the dense amalgam of the CAT(0) boundaries $\partial\Delta_i$, i.e. $\partial\Delta\cong\widetilde\sqcup(\partial\Delta_{1},\dots,\partial \Delta_{k})$.} 
\item{(4)} {\it Suppose that all vertex groups $G_{v_i}$ of $\cal G$ are systolic
(in the sense of simplicial nonpositive curvature as introduced by T. Januszkiewicz and the author in [JS]), and for each $v_i$ let $\Sigma_{i}$ be a systolic simplicial complex on which $G_{v_i}$ acts geometrically. Then there is a systolic complex $\Sigma$ on which $G$ acts geometrically, and such that the systolic boundary $\partial\Sigma$ (as introduced by D. Osajda and P. Przytycki in [OP])  is homeomorphic to the dense amalgam of the systolic boundaries $\partial\Sigma_i$, i.e. $\partial \Sigma\cong\widetilde\sqcup(\partial\Sigma_{1},\dots,\partial\Sigma_{k})$.} 

\medskip
Note that parts (2), (3) and (4) of Theorem 0.3 do not follow automatically from part (1), as a group $G$ may have in general many pairwise
non-homeomorphic $E{\cal Z}$-boundaries.
The $E{\cal Z}$-boundaries as in the assertion of part (1) 
have been constructed recently by C. Tirel [Ti] (the case of the free product) and by A. Martin [Ma] (the general case). We provide identifications of these boundaries  with the appropriate dense amalgams by referring to the characterization given in Theorem 0.2.
Part (2) of Theorem 0.3 strengthens an earlier result of A. Martin and the author [MS] (saying that, as a topological space, $\partial G$ depends
uniquely on the topology of the boundaries $\partial G_{v_i}$);
the strengthening concerns recognizing $\partial G$ as the appropriate dense amalgam. 

\smallskip
In Section 7 we  present
a more specific result concerning boundaries of Coxeter groups. Recall that to any Coxeter system $(W,S)$ there is associated
a CAT(0) polyhedral complex called the Coxeter-Davis complex
(see Chapters 7 and 12 in [Da], where this complex is denoted
by $\Sigma(W,S)$).
The CAT(0) boundary of this complex, denoted $\partial(W,S)$, is what is
shortly called the boundary of the Coxeter group $W$
(though it actually depends also on $S$).

\medskip\noindent
{\bf 0.4 Theorem.}
{\it Let $(W,S)$ be a Coxeter system, and suppose that $W$ has infinitely
many ends, and is not virtually free. Let $(W_1,S_1),\dots,(W_k,S_k)$
be the Coxeter systems corresponding to all maximal 1-ended
special subgroups of $W$. Then $k\ge1$, and
$$
\partial(W,S)\cong\widetilde\sqcup
\big(  \partial(W_1,S_1),\dots,\partial(W_k,S_k) \big).
$$}
\indent
The reader is advised to look also at the statement of Proposition 7.3.2
in the text. This proposition is the main step in the proof of Theorem 0.4,
but it also nicely complements the picture of appearance of 
dense amalgams as boundaries of Coxeter groups.

As it will be explained in Section 7, assumptions of Theorem 0.4
can be easily verified in terms of the Coxeter matrix of the system $(W,S)$.
Similarly, maximal 1-ended special subgroups of $W$ are easy to list
in terms of the same data. Note that Theorem 0.4 concerns all 
Coxeter systems except those for which the corresponding group
$W$ is either finite, or 2-ended, or 1-ended, or virtually free.
Thus, up to understanding the boundaries in 1-ended cases,
the theorem
presents a complete insight into the topology of boundaries
of Coxeter groups.

\medskip
The author thanks Krzysztof Omiljanowski
for helpful discussions.

\vfill\break

\bigskip

\noindent
{\bf 1. The dense amalgam.}

\medskip
In the major initial part of this section, given a nonempty 
compact metric space $X$, we construct out of it the unique
(up to homeomorphism) compact metric space $Y=\widetilde\sqcup(X)$,
called {\it the dense amalgam} of (copies of) $X$, and we show that
it satisfies conditions (a1)-(a5) of the introduction (for parameter $k=1$).
The construction of the space $\widetilde\sqcup(X)$ is rather involved
and requires a lot of auxilliary terminology and preparatory observations.

In the short final part of the section we extend the construction
to describe the dense amalgam $\widetilde\sqcup(X_1,\dots,X_k)$ of a finite
collection of compact metric spaces. 

\bigskip\noindent
{\it The peripheral extension $\ddot{X}$ of $X$.}

\medskip
Denote by $P$ the infinite countable discrete topological space.
Given a compact metric space $X$, its {\it peripheral extension}
is a compact metric space $K$ which contains $P$ as an open dense 
subspace such that $K\setminus P$ is homeomorphic to $X$.
In other words, $K$ is a metric compactification of $P$ with
the remainder $X$. Points of $P$ are called {\it the peripheral points}
of the extension $K$.

\medskip\noindent
{\bf Example.} Let $G$ be an infinite word hyperbolic group, 
and let $\partial G$ be its Gromov boundary. Then $\overline{G}=G\sqcup\partial G$,
equipped with the Gromov boundary compactification topology,
is a peripheral extension of the boundary $\partial G$. Its peripheral
points are precisely the elements of $G$.

\medskip
We record the following rather easy observations.

\medskip\noindent
{\bf 1.1 Lemma.}

\item{(1)} {\it Any nonempty compact metric space $X$ admits a peripheral extension.}
\item{(2)} {\it Any two peripheral extensions of a nonempty compact metric space $X$
are homeomorphic {\rm rel} $X$ (i.e. via a homeomorphism that is identical
on $X$).}
\item{(3)} {\it Given a peripheral extension $K$ of $X$, the group
of homeomorphisms of $K$ identical on $X$ acts transitively on the set $P$
of peripheral points of $K$.}

\medskip
In view of parts (1) and (2) of Lemma 1.1, a space $K$ as above exists and is uniquely
determined by $X$, so we denote it by $\ddot{X}$ and call {\it the peripheral
extension of $X$}.

\bigskip\noindent
{\it Complete tree systems of peripheral extensions of $X$.}

\medskip
Denote by $T$ the unique up to isomorphism countable tree
of infinite valence at every vertex. Let $V_T$ be the vertex set of $T$,
and $O_T$ the set of all oriented edges of $T$.
For each $t\in V_T$, we denote by $N_t$ the set of all oriented edges
of $T$ with initial vertex $t$.

A {\it complete tree system of peripheral extensions of $X$}
is a tuple $\Theta=(\{X_t\},\{b_t\})$ such that
to each $t\in V_T$ there is associated
\item{$\bullet$}  a space $X_t$
homeomorphic to $X$, equipped with its peripheral extension
$\ddot{X}_t$, and with the set $P_t$ of peripheral points of this extension;
\item{$\bullet$} a bijective map $b_t:N_t\to P_t$.

\noindent
Given two complete tree systems $\Theta=(\{X_t\},\{b_t\})$ and
$\Theta'=(\{X_t'\},\{b_t'\})$ of peripheral extensions of $X$,
an {\it isomorphism} between them is a tuple $F=(\lambda,\{ f_t \})$
such that:
\item{(I1)} $\lambda:T\to T$ is an automorphism;
\item{(I2)} for each $t\in V_T$ the map $f_t:\ddot{X}_t\to\ddot{X}_{\lambda(t)}'$ is a homeomorphism of peripheral extensions, i.e. it maps $X_t$ on $X_{\lambda(t)}'$ (and thus also $P_t$ onto $P_{\lambda(t)}'$);
\item{(I3)} for each $t\in V_T$ and any $e\in N_t$ the following commutation rule holds: $$b'_{\lambda(t)}(\lambda(e))=f_t(b_t(e)).$$

\medskip\noindent
{\bf 1.2 Lemma.}
{\it Any two complete tree systems of peripheral extensions of $X$
are isomorphic.}

\medskip\noindent
{\bf Proof:}
Let $\Theta=(\{X_t\},\{b_t\})$ and $\Theta'=(\{X_t'\},\{b_t'\})$ 
be two complete tree systems of peripheral extensions of $X$. 
Order the vertices of $V_T$ into a sequence $t_0,t_1,\dots$ 
so that for any natural $k$ the subtree $T_k$ of $T$ spanned 
on the vertices $t_1,\dots,t_k$ contains no other vertices of $V_T$. 
We construct an isomorphism $\lambda:T\to T$ and
homeomorphisms $f_t:\ddot{X}_t\to \ddot{X}'_{\lambda(t)}$ successively, 
at vertices $t=t_k$, as follows.
For each subtree $T_k$, denote by $T_k^+$ the subtree of $T$
spanned on $T_k$ and all vertices adjacent to the vertices of $T_k$.
Choose any $t_0'\in V_T$ and any homeomorphism  
$f_{t_0}:\ddot{X}_{t_0}\to\ddot{X}'_{t_0'}$ of peripheral extensions
(which exists by Lemma 1.1(2)). Denote by $T_0'$ 
the subtree of $T$ 
reduced to the vertex $t_0'$.  
Consider the bijection 
$$
(b'_{t_0'})^{-1}f_{t_0}b_{t_0}:N_{t_0}\to N_{t_0'}
$$
and denote by $\lambda_0:T_0^+\to (T_0')^+$ the isomorphism induced 
by the assignment $t_0\to t_0'$ and by the above bijection.

Now, suppose that we have already chosen the following data:
\item{(1)} vertices $t_0',\dots,t_k'$ in $V_{T}$ 
such that the subtree $T_k'$ of $T$ spanned on these vertices 
contains no other vertices of $V_{T}$, and the assignements 
$t_i\to t_i'$ yield an isomorphism $\mu_k:T_k\to T_k'$;
\item{(2)} an isomorphism $\lambda_k:T^+_k\to (T'_k)^+$ which
extends $\mu_k$;
\item{(3)} 
for $i=0,1,\dots,k$, homeomorphisms $f_{t_i}:\ddot{X}_{t_i}\to\ddot{X}'_{t_i'}$ 
of peripheral extensions such that the bijections 
$(b'_{t_i'})^{-1}f_{t_i}b_{t_i}:N_{t_i}\to N_{t_i'}$ 
are consistent with $\lambda_k$ (i.e. coincide with the
appropriate restrictions of the map induced by $\lambda_k$
between the sets of oriented edges of $T_k^+$ and $(T_k')^+$).

\noindent
Consider the vertex $t_{k+1}$, and 
let $j\in\{ 0,1,\dots,k \}$ be the index for which $t_j$
is the unique vertex of $T_k$ adjacent to $t_{k+1}$.
Put $t'_{k+1}:=\lambda_k(t_{k+1})$. 
Applying Lemma 1.1(3), choose any homeomorphism of peripheral extensions 
$f_{t_{k+1}}:\ddot{X}_{t_{k+1}}\to\ddot{X}'_{t'_{k+1}}$ such that 
$f_{t_{k+1}}(b_{t_{k+1}}([t_{k+1},t_j]))=b'_{t'_{k+1}}([t'_{k+1},t'_j])$. 
Denote by $\lambda_{k+1}:T_{k+1}^+\to (T'_{k+1})^+$ the isomorphism 
induced by $\lambda_k$ and the bijection 
$(b'_{t_{k+1}'})^{-1}f_{t_{k+1}}b_{t_{k+1}}:N_{t_{k+1}}\to N_{t_{k+1}'}$.

Iterating the above described step of the construction, we get an isomorphism 
$\lambda=\cup_k\lambda_k:T\to T$ and a family of homeomorphisms of peripheral extensions $f_t:\ddot{X}_t\to\ddot{X}'_{\lambda(t)}$ such that for each $t\in V_T$ the map
$(b'_t)^{-1}f_tb_t:N_t\to N_{t'}$ is consistent with $\lambda$ 
(i.e. coincides with the restriction of $\lambda$ to $N_t$). Since the latter clearly implies the commutativity condition (I3), we get that $F=(\lambda,\{ f_t \}):\Theta\to\Theta'$ is an isomorphism, which completes the proof.

\vfill\break

\bigskip\noindent
{\it The dense amalgam of (copies of) $X$.}

\medskip
We first describe an auxilliary compact metrisable space,
uniquely determined by $X$ up to homeomorphism, which intuitively
is the infinitely iterated and appropriately completed wedge of copies 
of $\ddot{X}$, in which the successively glued copies have rapidly decreasing size; 
wedge gluings are performed
at all peripheral points in all copies of $\ddot{X}$ so that exactly two 
copies meet at each gluing point.

More precisely, let $\Theta=(\{X_t\},\{b_t\})$ be a complete tree system
of peripheral extensions of $X$. 
For any finite subtree $F$ of $T$ define the {\it partial wedge of $\Theta$
for $F$}, 
as the quotient topological space
$$
\vee_F\Theta:=\sqcup_{t\in V_F}\ddot{X}_t\slash\sim,
$$
where $\sim$ is the equivalence relation induced by the equivalences
$$
b_t([t,s])\sim b_s([s,t])
$$
for all oriented edges $[t,s]$ of $F$.

For any pair of finite subtrees of $T$ such that $F'\subset F$
view $\vee_{F'}\Theta$ canonically as a subset of $\vee_F\Theta$, and
consider the retraction map $\rho_{F,F'}:\vee_F\Theta\to\vee_{F'}\Theta$
determined by the following. For any vertex $s\in V_F\setminus V_{F'}$
and for any $x\in\ddot{X}_s$, viewing $\ddot{X}_s$ canonically as a subset
of $\vee_F\Theta$, we put 
$$\rho_{F,F'}(x)=b_t([t,t']),$$
where $[t,t']$ is the last oriented edge on the shortest path in $T$
connecting $s$ with $F'$. Clearly, the retraction map $\rho_{F,F'}$ is continuous. 
Moreover, it is easy to check that for any finite subtrees
$F''\subset F'\subset F$ of $T$ we have
$\rho_{F',F''}\circ\rho_{F,F'}=\rho_{F,F''}$.

\medskip\noindent
{\bf 1.3 Definition.} 
\item{(1)}
{\it 
The wedge inverse system associated to $\Theta$} is the system over the poset of finite subtrees of $T$ given by
$$
{\cal S}^\vee_\Theta=(\{ \vee_F\Theta: F\subset T\hbox{ is a finite subtree} \}, 
\{ \rho_{F,F'}:F'\subset F\subset T \}).
$$
\item{(2)} 
{\it 
The wedge of $\Theta$} is the inverse limit of the system ${\cal S}^\vee_\Theta$, $$\vee\Theta:=\lim_{\longleftarrow}{\cal S}^\vee_\Theta.$$

\medskip
Since all partial wedges $\vee_F\Theta$ are easily seen to be compact 
metrisable, the same is true for their inverse limit $\vee\Theta$.

\medskip
Before getting further, we need to distinguish the subset $P_\Theta$ 
in $\vee\Theta$ consisting of the "gluing points" of the wedge.
More precisely, for any oriented edge $e=[t,t']\in O_T$ the point
$b_t([t,t'])\in\ddot{X}_t$, viewed as a point of $\vee\Theta$, coincides
with the point $b_{t'}([t',t])$, and we denote the corresponding point
of $\vee\Theta$ by $p_{|e|}$ (to emphasise the fact that it is induced
by the underlying non-oriented edge $|e|$). We then put
$$
P_\Theta:=\{ p_{|e|}:e\in O_T \}.
$$

\medskip\noindent
{\bf 1.4 Lemma.}
{\it Each point of the subset $P_\Theta$ is isolated in the
space $\vee\Theta$. In particular, $P_\Theta$ is an open subset
in $\vee\Theta$, and thus its complement $\vee\Theta\setminus P_\Theta$
is a compact metrisable space.}

\medskip\noindent
{\bf Proof:} Let $p=p_{|e|}$ be any point of $P_\Theta$.
Viewing $|e|$ as a subtree of $T$, we clearly have 
$p_{|e|}\in\vee_{|e|}\Theta\subset\vee\Theta$. Moreover,
if we denote by $\rho_{|e|}:\vee\Theta\to\vee_{|e|}\Theta$ the map
canonically associated to the inverse limit, it is not hard to see
that $\rho_{|e|}^{-1}(p_{|e|})=p_{|e|}$. Since $p_{|e|}$ is 
isolated in $\vee_{|e|}\Theta$, its singleton is an open subset
in $\vee_{|e|}\Theta$, and thus the same is true in $\vee\Theta$,
which completes the proof.

\medskip
Note that, it follows easily from the above description of $\vee\Theta$
that if $\Theta$ and $\Theta'$ are two isomorphic complete tree systems
of peripheral extensions of $X$, then the pairs of spaces $(\vee\Theta,P_\Theta)$
and $(\vee\Theta',P_{\Theta'})$ are homeomorphic.
This and Lemma 1.2 then justify the following.

\medskip\noindent
{\bf 1.5 Definition.}
The {\it dense amalgam of (copies of) $X$}, denoted $\widetilde\sqcup(X)$, 
is the topological space
$\vee\Theta\setminus P_\Theta$, where $\Theta$ is any complete
tree system of peripheral extensions for $X$.

\bigskip\noindent
{\it A more explicit description of the wedge $\vee\Theta$ and its subspace
$\vee\Theta\setminus P_\Theta$.}

\medskip
Given a complete tree system $\Theta=(\{ \ddot{X}_t \},\{ b_t \})$ 
of peripheral extensions
of $X$, consider the equivalence relation on the disjoint union
$\sqcup_{t\in V_T}\ddot{X}_t$ induced by the equivalences
$b_t([t,s])\sim b_s([s,t])$
for all oriented edges $[t,s]\in O_T$.
Denote the set of equivalence classes of this relation by $\#\Theta$.
Let $\partial T$ be the set of ends of the tree $T$,
i.e. the set of equivalence classes for the relation 
on the set of infinite rays in $T$
provided by coincidence of two rays except possibly at
some finite initial part in each of them.
Since the inverse system ${\cal S}^\vee_\Theta$ consists of natural
retractions of bigger partial wedges on the smaller ones,
one easily identifies the inverse limit $\vee\Theta$, set theoretically,
with the disjoint union $\#\Theta\sqcup\partial T$.

We now describe the topology of the inverse limit $\vee\Theta$
as topology on the set $\#\Theta\sqcup\partial T$.
For any finite subtree $F$ of $T$ consider the map
$\rho_F:\vee\Theta\to\vee_F\Theta$ canonically associated to the
inverse limit. 
Under identification of $\vee\Theta$ with $\#\Theta\sqcup\partial T$,
this map is easily seen to have the following form.
If $s\in V_F$ and $x\in\ddot{X}_s\subset\#\Theta$, then
$\rho_F(x)=x\in\ddot{X}_s\subset\vee_F\Theta$.
If $s\in V_T\setminus V_F$, let $[t_s,t'_s]$ be the first oriented edge
on the unique minimal path connecting a vertex of $F$ to $s$; then for any 
$x\in \ddot{X}_s$ we have $\rho_F(x)=b_{t_s}([t_s,t'_s])\in\ddot{X}_{t_s}\subset\vee_F\Theta$.
Finally, if $z\in\partial T$, let $[t_z,t'_z]$ be the first oriented edge
on the unique minimal ray in $T$ representing $z$ and starting at a vertex of $F$;
then $\rho_F(z)=b_{t_z}([t_z,t'_z])\in\ddot{X}_{t_z}\subset\vee_F\Theta$.

By definition of the inverse limit, the family
$$
\{ \rho_F^{-1}(U):F \hbox{ is a finite subtree of }T \hbox{ and }
U \hbox{ is an open subset of }\vee_F\Theta \}
$$
is a subbasis for the topology in $\vee\Theta$.
It follows from the above description of $\rho_F$ that
any subset $\rho_F^{-1}(U)$ from this subbasis, viewed as a subset of $\#\Theta\sqcup\partial T$,
can be described as follows. Identify $\vee_F\Theta$
and all the sets $\ddot{X}_t$ canonically as the subsets in $\#\Theta$.
Under notation as in the previous paragraph,
put 
$$
\#_U\Theta:=U\cup\bigcup\{ \ddot{X}_s:s\in V_T\setminus V_F \hbox{ and }
b_{t_s}([t_s,t'_s])\in U \}\subset\#\Theta.
$$
Furthermore, put
$$
\partial_U T:=\{ z\in\partial T:b_{t_z}([t_z,t'_z])\in U \}.
$$
Then $\rho_F^{-1}(U)=\#_U\Theta\cup\partial_U T$.

\medskip
Using the above description of the sets $\rho_F^{-1}(U)$,
we now indicate a convenient basis of the topology in 
$\vee\Theta=\#\Theta\sqcup\partial T$. For any vertex $t\in V_T$,
viewing it as a subtree of $T$, we denote by
$\rho_t:\vee\Theta\to\vee_t\Theta=\ddot{X}_t$ the map canonically
associated to the inverse limit.

\medskip\noindent
{\bf 1.6 Lemma.}
{\it The family
$$
{\cal B}=
\{ \rho_t^{-1}(U):t\in V_T \hbox{ and }U \hbox{ is an open subset of }\ddot{X}_t \} \cup \{ \{p\}:p\in P_\Theta \}
$$
is a basis of the topology in $\vee\Theta$.}

\medskip\noindent
{\bf Proof:}
We will first show that the family $\cal B$ satisfies the axioms 
of a basis of topology. 
Since $\cal B$ is obviously a covering of $\vee\Theta$,
it remains to check that the intersection $B\cap B'$ of any two sets from $\cal B$
is the union of some sets from $\cal B$. This is obvious if $B$ or $B'$
is a singleton from $P_\Theta$. Thus, we need to study the case
when $B=\rho_t^{-1}(U)$ and $B'=\rho_s^{-1}(U')$,
where $U,U'$ are some open subsets in $\ddot{X}_t$ and $\ddot{X}_s$, respectively. 

If $t=s$, we get $B\cap B'=\rho_t^{-1}(U\cap U')$, which trivially yielkds our assertion. If $t\ne s$, 
let $F$ be the subtree of $T$ spanned on $t$ and $s$
(which is obviously finite), and let $O_F$ be the set of oriented
edges in $F$.
Put $U_t:=U\setminus b_t(N_t\cap O_F)$, 
$U_s:=U\setminus b_s(N_t\cap O_F)$, and for each 
$a\in V_F\setminus\{ t,s \}$ put
$U_a:=\ddot{X}_a\setminus b_a(N_a\cap O_F)$.
Observe that for any $a\in V_F$ the set $N_a\cap O_F$ is finite.
Consequently, for any $a\in V_F$ the set $b_a(N_a\cap O_F)$ is closed,
and hence $U_a$ is open in the corresponding space $\ddot{X}_a$.
Furthermore, define a subset $A\subset V_F$ by the following rules:
\item{$\bullet$} $t$ belongs to $A$ if $\ddot{X}_t\subset B'$,
\item{$\bullet$} $s$ belongs to $A$ if $\ddot{X}_s\subset B$,
\item{$\bullet$} a vertex $a\in V_F\setminus\{ t,s \}$ belongs to $A$ if $\ddot{X}_a\subset B\cap B'$.

\noindent
It is not hard to observe that 
$$
B\cap B'=(P_\Theta\cap B\cap B')\cup\bigcup_{a\in A}\rho_a^{-1}(U_a),
$$
which also yields our assertion. Thus 
$\cal B$ satisfies the axioms 
of a basis of topology. 

Now we need to show that the topology ${\cal T}_{\cal B}$
induced by $\cal B$ coincides with the original topology $\cal T$
in $\vee\Theta=\#\Theta\sqcup\partial T$. Since, in view of Lemma 1.4
we have ${\cal B}\subset{\cal T}$, it follows that 
${\cal T}_{\cal B}\subset{\cal T}$.
To prove the converse inclusion, it is enough to show
that any set of form $\rho_F^{-1}(U)$, where $F$ is any finite subtree of $T$, is the union of some
elements of $\cal B$.
To do this, for each $a\in V_F$ put 
$U_a:=(U\cap\ddot{X}_a)\setminus b_a(N_a\cap O_F)$.
Note that, by the argument as before, this is an open subset
of $\ddot{X}_a$. Observe that we have
$$
\rho_F^{-1}(U)=(P_\Theta\cap\rho_F^{-1}(U))\cup
\bigcup_{a\in V_F}\rho_a^{-1}(U_a),
$$
which completes the proof.

\medskip
We now pass to the subspace $\vee\Theta\setminus P_\Theta$.
Consider the family of its subsets
$$
{\cal B}_0:=\{ W\setminus P_\Theta:W\in{\cal B} \}=
\{ \rho_t^{-1}(U)\setminus P_\Theta:t\in V_T \hbox{ and }U 
\hbox{ is an open subset in }\ddot{X}_t \}.
$$

From Lemma 1.6 we immediately get the following.

\medskip\noindent
{\bf 1.7 Corollary.}
{\it ${\cal B}_0$ is a basis of the topology in $\vee\Theta\setminus P_\Theta$.}

\medskip
Note that, under identification $\vee\Theta=\#\Theta\sqcup\partial T$,
the subspace $\vee\Theta\setminus P_\Theta$ is identified with the
subset $(\bigsqcup_{t\in V_T}X_t)\sqcup\partial T$. By what was said above,
we have the following description of any set 
$\rho_t^{-1}(U)\setminus P_\Theta\in{\cal B}_0$
as a subset of $(\bigsqcup_{t\in V_T}X_t)\sqcup\partial T$.
For $s\in V_T\setminus\{t\}$ let $[t,t_s]$ be the first oriented edge
on the path in $T$ from $t$ to $s$. Similarly, for any $z\in\partial T$
let $[t,t_z]$ be the fist oriented edge on the unique ray in $T$
started at $t$ and representing $z$. Recalling that $U$ is an open subset
of $\ddot{X}_t$, define the subset $D(t,U)\subset (\bigsqcup_{t\in V_T}X_t)\sqcup\partial T$
as
$$
D(t,U):=(U\cap X_t)\sqcup\bigsqcup\{ X_s:s\ne t \hbox{ and } b_t([t,t_s])\in U \}
\sqcup\{ z\in\partial T:b_t([t,t_z])\in U \}.
$$
We then have $\rho_t^{-1}(U)\setminus P_\Theta=D(t,U)$.

As immediate restatement of Corollary 1.7 we get the following.

\medskip\noindent
{\bf 1.8 Proposition.}
{\it The family 
$$
{\cal D}=\{ D(t,U):t\in V_T \hbox{ and }U \hbox{ is an open subset of }\ddot{X}_t \}
$$
is a basis for the topology in the space $\vee\Theta\setminus P_\Theta$,
under its identification with 
$(\bigsqcup_{t\in V_T}X_t)\sqcup\partial T$.}

\medskip
To conclude the explicit description of the space $\vee\Theta\setminus P_\Theta$
(and thus also the amalgam $\widetilde\sqcup(X)$),
we provide in the lemma below 
some bases of open neighbourhoods for all points in this space.

\medskip\noindent
{\bf 1.9 Lemma.}
{\it For the canonical identification of the space $\vee\Theta\setminus P_\Theta$
with $(\bigsqcup_{t\in V_T}X_t)\sqcup\partial T$, we have:}
\item{(1)} {\it if $x\in X_t\subset (\bigsqcup_{t\in V_T}X_t)\sqcup\partial T$, 
then the family of sets $D(t,U)$, where $U$ runs through any basis of
open neighbourhoods of $x$ in $\ddot{X}_t$, is a basis of open neighbourhoods
of $x$ in $(\bigsqcup_{t\in V_T}X_t)\sqcup\partial T$;}
\item{(2)} {\it if $z\in\partial T\subset (\bigsqcup_{t\in V_T}X_t)\sqcup\partial T$,
then for any ray $[t_0,t_1,\dots]$ in $T$ representing $z$ the family
$$
\{ D(t_i,\ddot{X}_{t_i}\setminus\{ b_{t_i}([t_i,t_{i-1}]) \}):i\ge1 \}
$$
is a basis of open neighbourhoods
of $z$ in $(\bigsqcup_{t\in V_T}X_t)\sqcup\partial T$.}

\medskip
We skip a straighforward proof of this lemma.

\bigskip\noindent
{\it The amalgam $\widetilde\sqcup(X)$ satisfies conditions (a1)-(a5).}

\medskip
As we have already noticed, the amalgam $\widetilde\sqcup(X)$
is a compact metrisable space. 
We now check that it satisfies conditions (a1)-(a5) listed in the introduction.
To do this, we will use the above discussed identification
of the space $\widetilde\sqcup(X)\cong\vee\Theta\setminus P_\Theta$ with the set
$(\bigsqcup_{t\in V_T}X_t)\sqcup\partial T$ equipped with topology
provided by the basis $\cal D$, as stated
in Proposition 1.8. As a family $\cal Y$ of subsets we take the family
$X_t:t\in V_T$. 

\smallskip
Note that each $X_t$ is an embedded
copy of $X$, as it coincides with the image of the canonical embedding of $X_t$
in the inverse limit $\vee\Theta$. Since the subsets in this family
are clearly pairwise disjoint, condition (a1) is fulfilled.

\smallskip
To check condition (a2) we need to verify that for any finite open
covering $\cal U$ of $(\bigsqcup_{t\in V_T}X_t)\sqcup\partial T$
and for almost every vertex $t\in V_T$ (i.e. for each $t\in V_T\setminus A$,
where $A$ is some finite subset of $V_T$) there is $U\in{\cal U}$
such that $X_t\subset U$. Obviously, without loss of generality,
we may assume that $\cal U$ consists of subsets from the basis $\cal D$.
More precisely, we may assume that there is a finite subset $A\subset V_T$
and a family $U_s:s\in A$ of open subsets in the corresponding spaces
$\ddot{X}_s$ such that ${\cal U}=\{ D(s,U_s):s\in A \}$.
But then it is easy to check that for each $t\in V_T\setminus A$
we have $X_t\subset D(s,U_s)$ for some $s\in A$, which verifies (a2).

\smallskip
Condition (a3) follows easily from the description
of bases of open neighbourhoods of points in 
$(\bigsqcup_{t\in V_T}X_t)\sqcup\partial T$,
as given in Lemma 1.9(1).
We skip this strightforward argument.
Similarly, condition (a4) follows
directly from Lemma 1.9(2).

\smallskip
To check condition (a5),  we introduce
a family of $\cal Y$-saturated open and closed subsets 
of  $(\bigsqcup_{t\in V_T}X_t)\sqcup\partial T$
that we call {\it half-spaces}.
For any edge $e=[t,t']\in V_T$ consider the subsets
$H_e^-:=D(t,\ddot{X}_t\setminus\{ b_t([t,t']) \})$ 
and $H_e^+:=D(t',\ddot{X}_{t'}\setminus\{ b_{t'}([t',t]) \})$,
and note that they are both open. Moreover, they form a partition
of the space $(\bigsqcup_{t\in V_T}X_t)\sqcup\partial T$,
and thus they are both open and closed. Finally, both these
subsets are easily seen to be $\cal Y$-saturated.
We will call them the {\it half-spaces} induced by the edge $e$.

Now, let $x,y$ be any two distinct points of $(\bigsqcup_{t\in V_T}X_t)\sqcup\partial T$ which do not belong to the same
subset of $\cal Y$.
First, consider the case when
$x\in X_t$ for some $t\in V_T$.
If $y\in X_s$ for some $s\ne t$, then for any oriented edge $e$
on the path connecting $t$ with $s$ we heve $A\subset H_e^-$
and $y\in H_e^+$. If $y\in\partial T$, then for any oriented edge
$e$ in the ray started at $t$ and representing $y$ we similarly
have $A\subset H_e^-$ and $y\in H_e^+$. This verifies  condition (a5)
in the considered case.
Since in the remaining case, when $x,y\in\partial T$,
we can also easily separate $x$ from $y$ 
by a half-space, condition (a5) follows.

\bigskip\noindent
{\it The dense amalgam $\widetilde\sqcup(X_1,\dots,X_k)$.}

\medskip
Given a finite collection $X_1,\dots,X_k$ of nonempty compact metric spaces,
put
$$
\widetilde\sqcup(X_1,\dots,X_k):=\widetilde\sqcup(X),
$$
where $X=X_1\sqcup\dots\sqcup X_k$ is the topological disjoint union.
Under identification 
$\widetilde\sqcup(X)=(\bigsqcup_{t\in V_T}X_t)\sqcup\partial T$,
for each $t\in V_T$ we have $X_t=X_{1,t}\sqcup\dots\sqcup X_{k,t}$,
where each $X_{i,t}$ is homeomorphic to the corresponding $X_i$.
For each $i\in\{ 1,\dots,k \}$  take ${\cal Y}_i:=\{ X_{i,t}:t\in V_T \}$,
and put ${\cal Y}:=\cup_i{\cal Y}_i$.
We check that so defined space and
the family of its subspaces satisfy conditions (a1)-(a5)
from the introduction.

The only condition which does not follow by an argument similar
as before is condition (a5), in the case of two points $x\in X_{i,t}$
and $y\in X_{j,t}$ for some $t\in V_T$ and some $j\ne i$.
Observe that we can partition the peripheral extension $\ddot{X}_t$
of $X_t=X_{1,t}\sqcup\dots\sqcup X_{k,t}$ into open and closed subsets $U,W$
such that $U\cap X_t=X_{i,t}$. The subsets $D(t,U)=\rho_t^{-1}(U)$
and $D(t,W)=\rho_t^{-1}(W)$ form then an open and closed partition of 
$\widetilde\sqcup(X_1,\dots,X_k)$, and since we obviously have
that $x\in X_{i,t}\subset D(t,U)$ and $y\in X_{j,t}\subset D(t,W)$,
the argument is completed.

\vfill\break

\bigskip

\noindent
{\bf 2. The characterization.}

\medskip
The aim of this section is to prove Theorem 0.2 of the introduction.
We start with introducing a useful terminology. 
Let $X_1,\dots,X_k$ be a collection of nonempty metric compacta,
for some $k\ge1$.
A compact metric space $Y$ is $(X_1,\dots,X_k)$-{\it regular} if it can be equipped
with a family $\cal Y$ of subspaces satisfying conditions (a1)-(a5)
from the introduction. Any family $\cal Y$ with these properties
is called an $(X_1,\dots,X_k)$-{\it regularizing} family for $Y$.
Theorem 0.2 may be then rephrased as follows: 
any $(X_1,\dots,X_k)$-regular compact metric space is homeomorphic
to the dense amalgam $\widetilde\sqcup(X_1,\dots,X_k)$.

In view of the definition of the dense amalgam 
$\widetilde\sqcup(X_1,\dots,X_k)$ for $k>1$, given at the end of Section 1,
Thorem 0.2 is a direct consequence of the following two results.

\medskip\noindent
{\bf 2.1 Proposition.}
{\it Given a nonempty compact metric space $X$, each $(X)$-regular
space $Y$ is homeomorphic to the dense amalgam $\widetilde\sqcup(X)$.}

\medskip\noindent
{\bf 2.2 Proposition.}
{\it Given any tuple $X_1,\dots,X_k$ of nonempty compact metric spaces,
each $(X_1,\dots,X_k)$-regular space $Y$ is also 
$(X_1\sqcup\dots\sqcup X_k)$-regular.}

\medskip
In the proofs of both propositions above we will use the following notation.
Given a nonempty subset $A$ in a metric space $Y$, and a real number
$\epsilon>0$, an {\it $\epsilon$-neighbourhood} of $A$ is the subset
$$
N_\epsilon(A):=\{ x\in Y:d_Y(A,x)<\epsilon \}.
$$

\noindent
The {\it diameter} of $A$ is the number 
$\hbox{diam}(A):=\sup\{ d_Y(x,y):x,y\in A \}$.

\medskip\noindent
{\bf Proof of Proposition 2.2:}
Let $\cal Y$ be an $(X_1,\dots,X_k)$-regularizing family for $Y$.
We will construct inductively
an $(X_1\sqcup\dots\sqcup X_k)$-regularizing family ${\cal W}=(W_n)_{n\ge1}$ for $Y$. Each  $W_n\in{\cal W}$ will have a form of the union
of some appropriately chosen subsets from $\cal Y$.

Order the elements of $\cal Y$ into
a sequence $(Y_n)_{n\ge1}$. Put $Z_{1,1}=Y_1$ and choose the subsets
$Z_{1,2},\dots,Z_{1,k}\in{\cal Y}$ such that:
\itemitem{(z1)} the family $Z_{1,i}:1\le i\le k$ consists of exactly
one set from each of the subfamilies ${\cal Y}_i$ of $\cal Y$;
\itemitem{(z2)} for each $2\le i\le k$ we have 
$Z_{1,i}\subset N_{{\rm diam}(Z_{1,1})}(Z_{1,1})$.

\noindent
Such a choice is possible since, by conditions (a2)-(a4), 
each family ${\cal Y}_i$ is infinite , null and dense in $Y$.
Put $W_1=Z_{1,1}\cup\dots\cup Z_{1,k}$. 

Having already constructed the subsets $W_1,\dots,W_{n-1}$ as
unions of some subfamilies of $\cal Y$, we construct the subset $W_n$
as follows. If $Y_n$ is not contained in $W_1\cup\dots\cup W_{k-1}$,
put $Z_{n,1}=Y_n$; otherwise, take as $Z_{n,1}$ any subset from $\cal Y$
not contained in $W_1\cup\dots\cup W_{k-1}$.
Choose $Z_{n,2},\dots,Z_{n,k}\in{\cal Y}$ not contained in 
$W_1\cup\dots\cup W_{k-1}$ and satisfying the analogons of
conditions (z1) and (z2) above, with $Z_{1,i}$'s replaced by $Z_{n,i}$'s.

We now check that $\cal W$ is an $(X_1\sqcup\dots\sqcup X_k)$-regularizing
family of subsets for $Y$, i.e. it satisfies the appropriate variant
of conditions (a1)-(a5).
Note that $\cal W$ obviously consists of subsets
which are embedded copies of $X_1\sqcup\dots\sqcup X_k$,
and each such copy is boundary in $Y$ (as finite union of 
closed boundary subsets).
Moreover, by condition (z2) for each $n$ we have 
${\rm diam}(W_n)\le 3{\rm diam}(Z_{n,1})$, and thus the family
$\cal W$ is null. Finally, it follows from the above description
that for each $n$ we have $Y_n\subset W_1\cup\dots\cup W_n$,
and so we have that $\cup{\cal W}=\cup{\cal Y}$.
In particular, the union of the family $\cal W$ is dense in $Y$.

It remains to show that the family $\cal W$ satisfies
condition (a5).
However, in order to ensure that this is true, we need to add some ingredient  to the construction presented above. To describe this ingredient, for each $n\ge1$ consider the number 
$$
z_n:=\max\{ \hbox{diam}(Z):Z\in{\cal Y}, Z\i\hskip-9pt\slash \cup_{j=1}^n W_j \}
$$
and note that, since the family $\cal Y$ is null, we have $\lim_n z_n=0$. Now, in the above inductive construction of the subspaces $W_n$, for each $n$ we additionally choose a finite partition ${\cal Q}_n$ of $Y$ into $\cal Y$-saturated closed and open subsets
$Q_1^n,\dots,Q_{m_n}^n$,
such that

\itemitem{(q1)} for each $n$ we have $m_n\ge n$, and 

\hskip1.2cm (a) for each $j\in\{1,\dots,n\}$ we have 
$W_j\i Q^n_j\i N_{1/n}(W_j)$, 

\hskip1.2cm (b) for each $j\in\{n+1,\dots,m_n  \}$
we have $\hbox{diam}(Q_j^n)<z_n+1/n$;

\itemitem{(q2)} ${\cal Q}^{n+1}$ is a refinement of ${\cal Q}^n$ for each $n\ge1$;

\itemitem{(q3)} for each $n$ the subset $W_{n+1}$ is contained in one of 
the sets $Q^n_j\in{\cal Q}^n$.

\noindent
More precisely, at each step of the construction, after choosing a subspace $W_n$ we choose a partition ${\cal Q}^n$ satisfying (q1) and (q2), and then we choose $W_{n+1}$ satisfying (q3). The possibility to choose ${\cal Q}^n$ satisfying (q1) and (q2) follows from condition (a5)  for the family $\cal Y$, due to the following.

\medskip\noindent
{\bf Claim 1.} 
{\sl If $A$ is either a subspace from $\cal Y$ or a point from 
the subset $Y^\bullet:=Y\setminus\cup{\cal Y}$, then $\forall\epsilon>0$ there is a closed and open $\cal Y$-saturated set $Q$ such that $A\i Q\i N_\epsilon(A)$.}

\medskip
We skip a straightforward proof of Claim 1, indicating only that it
uses the fact that $\cal Y$-saturated closed and open subsets of $Y$ are closed
under finite intersections and finite unions. Once we have chosen ${\cal Q}^n$,
in the description of $W_{n+1}$ as above we additionally require that all the sets 
$Z_{n+1,i}:2\le i\le k$ are contained in the same $Q_j^n$ as the set 
$Z_{n+1,1}$, which guaranties (q3).

Observe that, by the above description, all closed and open subsets $Q^n_j$ appearing in any of the partitions ${\cal Q}^n$ are ${\cal W}$-saturated. Thus,
we may use them as separating sets justifying condition (a5).
Namely, if $x\in W_n$ for some $n\ge1$,  then $x$ can be separated from a point $y\notin W_n$ by a subset $Q^m_n$, for sufficiently large $m$, due to condition (q1)(a). If $x\in Y^\bullet$,  for each $n$ consider this $j_n$ for which $x\in Q^n_{j_n}$. 
We will need the following.

\medskip\noindent
{\bf Claim 2.}
{\sl $\hbox{diam}(Q_{j_n}^n)\to0$.}

\medskip
To prove Claim 2, consider first the case when for each $n$ we have
$j_n\le n$. In this case, by condition (q1)(a), we have 
$x\in N_{1/n}(W_{j_n})$ for all $n$. From this it is not hard to deduce
that for a fixed $j$ we have $j_n=j$ only for finitely many $n$,
and hence $j_n\to\infty$. Since $Q^n_{j_n}\subset N_{1/n}(W_{j_n})$
and so $\hbox{diam}(Q^n_{j_n})\le\hbox{diam}(W_{j_n})+{2\over n}$,
we get that $\hbox{diam}(Q_{j_n}^n)\to0$ by the fact that the family
$\cal W$ is null.

Now, consider the case when for each $n$ we have $j_n>n$.
It follows that $j_n\to\infty$. By condition (q1)(b),
we have $\hbox{diam}(Q_{j_n}^n)<z_n+{1\over n}$,
and hence $\hbox{diam}(Q_{j_n}^n)\to0$ in this case too.
The general case easily follows from the two just considered cases,
hence Claim 2.

\smallskip
By Claim 2, $x$ can be separated from any other point $y\in Y$
by a set $Q_{j_n}^n$, for sufficiently large $n$.
This completes the proof of Proposition 2.2.

\bigskip
The proof of Proposition 2.1 requires more terminology and auxilliary
results, which we provide in four preparatory subsections below.
The proof itself appears at the end of the section.

In all the remaining part of this section we work under notation and assumptions of Proposition 2.1. It means that $X$ is a nonempty metric compactum, $Y$ is  an $(X)$-regular space, and $\cal Y$ is
an $(X)$-regularizing family for $Y$. We fix a metric $d_Y$ in $Y$.
We also often refer to the subset $Y^\bullet=Y\setminus\cup{\cal Y}$.

\bigskip\noindent
{\it 2.A Cantor space $C$ and the related space $C_0$.}

\medskip
Recall that the {\it Cantor space} is a metric compactum $C$ determined uniquely by the following properties:
\itemitem{(c1)} $C$ is zero-dimensional, i.e. every point of $C$ is a connected component of $C$ (it can be separated from any other point by a closed and open subset of $C$); 
\itemitem{(c2)} $C$ has no isolated points, i.e. every point of $C$ is an accumulation point.

\medskip\noindent
{\bf 2.A.1 Lemma.} 
{\sl The quotient space $Y/{\cal Y}$ is homeomorphic to the Cantor space $C$.}

\medskip\noindent
{\bf Proof:} Since $\cal Y$ is a null decomposition of $Y$, it follows from [Dav, Proposition 2 on p. 13] that $Y/{\cal Y}$ is a metric compactum. We need to check conditions (c1) and (c2). 
Condition (c2) follows easily from condition (a3) for $\cal Y$, and condition (c1) is a consequence of condition (a5), hence the lemma.

\medskip
We now recall or provide few properties of the Cantor space and its subspaces that will be useful later in this section. Denote by $C_0$ the space obtained by deleting any single point from the Cantor space $C$. The following two results are well known.

\medskip\noindent
{\bf 2.A.2 Proposition.}
{\sl A locally compact metric space is homeomorphic to $C_0$ iff it is zero-dimensional, noncompact and has no isolated points.}

\medskip\noindent
{\bf 2.A.3 Proposition.}
{\sl Any noncompact open subset of $C$ is homeomorphic to $C_0$. In particular, the complement $C\setminus Z$ of any closed boundary subset $Z\i C$ is homeomorphic to $C_0$.}

\medskip
We will also need the following technical result.

\medskip\noindent
{\bf 2.A.4 Lemma.}
{\sl Let $\{ p_\lambda:\lambda\in\Lambda \}$ be an infinite (in fact, countable) discrete subset of the space $C_0$, and let $U_\lambda:\lambda\in\Lambda$ be a covering of $C_0$ by open subsets with compact closures in $C_0$ such that $p_\lambda\in U_\lambda$ for each $\lambda\in\Lambda$. Then there is a partition of $C_0$ into subsets $K_\lambda:\lambda\in\Lambda$ which are compact, open, and such that $p_\lambda\in K_\lambda\i U_\lambda$ for each $\lambda\in\Lambda$.
Moreover, the subsets $K_\lambda$ are all homeomorphic to the Cantor space.}

\medskip\noindent
{\bf Proof:}
The first assertion is a fairly straightforward consequence of the fact
that each point of $C_0$ has a basis of open neighbourhoods which
are also compact. The second assertion follows from the fact that any nonempty closed and open subset of the Cantor space $C$ is
homeomorphic to $C$.

\vfill\break

\bigskip\noindent
{\it 2.B Sequences of subspaces convergent to points.}

\medskip

Since the family $\cal Y$ is null, given any infinite sequence $(Z_n)$ of pairwise distinct subspaces from $\cal Y$, we have $\lim_{n\to\infty}\hbox{diam}(Z_n)=0$. This justifies the following. Given a sequence $(Z_n)$ as above, we say that a point $p\in Y$ is the limit of this sequence, $\lim_{n\to\infty}Z_n=p$, if for some (and hence any) selection of points $p_n\in Z_n$ we have $\lim_{n\to\infty}p_n=p$.
In such a situation we also say that the sequence $(Z_n)$ is convergent.

\medskip\noindent
{\bf 2.B.1 Fact.}
{\sl Each point $p\in Y$ can be expressed as $p=\lim_{n\to\infty}Z_n$ for some sequence $(Z_n)$ as above.}

\medskip\noindent
{\bf Proof:} If $p\in Y^\bullet$, the assertion follows directly from condition (a4) for $\cal Y$ (and from compactness of the subspaces in $\cal Y$). If $p\in Z\in{\cal Y}$ then, by condition (a3) applied to $Z$, $p$ is either the limit as required, or the limit of some sequence $(p_n)$ of points from the subset $Y^\bullet$. In the latter case, since each $p_n$ is the limit as required, the same holds for $p=\lim_n p_n$, which completes the proof.

\medskip
We present two more technical results concerning convergent sequences of subspaces from $\cal Y$. We skip a straightforward proof of the first of these two results.

\medskip\noindent
{\bf 2.B.2 Lemma.}
{\sl Let $(Z_n)_{n\ge1}$ be a sequence of subspaces in $Y$ satisfying the following conditions:}
\item{(1)} {\sl $Z_1$ is arbitrary;}
\item{(2)} {\sl for each $n\ge1$ we have $Z_{n+1}\i N_{{\rm diam}(Z_n)}(Z_n)$ and $\hbox{\rm diam}(Z_{n+1})<{1\over2}\hbox{\rm diam}(Z_n)$.}

\noindent
{\sl Then $(Z_n)$ is convergent, and if $p=\lim_n Z_n$ then}
$$
\{ p \}\cup\bigcup_{n\ge1}Z_n\i N_{2{\rm diam}(Z_1)}(Z_1).
$$

\medskip\noindent
{\bf 2.B.3 Lemma.}
{\sl Given an ordering of the family ${\cal Y}$  into a sequence $\{ Y_n:n\in N \}$, let $(Z_n)_{n\ge0}$ be a sequence of distinct subspaces from $\cal Y$ satisfying the following conditions:}
\item{(0)} $Z_0$ is arbitrary;
\item{(1)} {\sl $Z_1\ne Y_1$, $\hbox{diam}(Z_1)<{1\over2}\hbox{diam}(Z_0)$ and $Z_1\i N_{{\rm diam}(Z_0)}(Z_0)$;}
\item{(2)} {\sl for each $n\ge1$ we have $Z_{n+1}\ne Y_{n+1}$,  $\hbox{\rm diam}(Z_{n+1})<{1\over2}\hbox{\rm diam}(Z_n)$ and $Z_{n+1}\i N_{d_n}(Z_n)$, where
$$
d_n:=\min\left( {\rm diam}(Z_n), {1\over3}d_Y(Y_n,Z_n), \dots, {1\over 3^{n}}d_Y(Y_1,Z_1) \right).
$$

\noindent
Then $(Z_n)$ is convergent, and the limit point $p=\lim_n Z_n$ belongs to $Y^\bullet$.}

\medskip\noindent
{\bf Proof:} Convergence follows from Lemma 2.B.2. Moreover, it is not hard to see that for the limit point $p$ we have $$d_Y(p,Y_n)>({1}-\sum_{i=1}^\infty{1\over 3^i})d_Y(Y_n,Z_n)={1\over2}d_Y(Y_n,Z_n)>0$$ for each $n\ge1$, and thus $p\notin\cup_{n=1}^\infty Y_n=\cup{\cal Y}$.

\bigskip\noindent
{\it 2.C Approximating families of subspaces.}

\smallskip
We will frequently use the following concept.

\medskip\noindent
{\bf 2.C.1 Definition.} Let $Z\in{\cal Y}$.
A subfamily ${\cal Y}_0\i{\cal Y}$ {\it approximates} $Z$ if:
\item{$\bullet$} $Z\notin{\cal Y}_0$;
\item{$\bullet$} $Z\i\overline{\cup{\cal Y}_0}$; 
\item{$\bullet$} $\lim_{W\in{\cal Y}_0}d_Y(W,Z)=0$ (equivalently, for any $\epsilon>0$ almost all $W\in{\cal Y}_0$ are contained in $N_\epsilon(Z)$).

\medskip
We make a record of few easily seen properties of approximating families.

\medskip\noindent
{\bf 2.C.2 Fact.}
\item{(1)} {\sl Any $Z\in{\cal Y}$ admits an approximating family ${\cal Y}_0\i{\cal Y}$.}
\item{(2)} {\sl Any approximating family ${\cal Y}_0$ for $Z\in{\cal Y}$ is {\it discrete} in the complement $Y\setminus Z$ in any of the following two equivalent senses:}
\itemitem{(a)} {\sl for each $W\in{\cal Y}_0$ there is $\delta>0$ such that the neighbourhood $N_\delta(W)$ (in $Y$) is disjoint with $Z$ and with all subspaces from ${\cal Y}_0\setminus\{W\}$;}
\itemitem{(b)} 
{\sl the subset $\{[W]:W\in{\cal Y}_0 \}\i Y/{\cal Y}$ is discrete in 
$(Y/{\cal Y})\setminus\{ [Z] \}$,
where for any
$W\in{\cal Y}$ we denote by $[W]$ the point in
the quotient $Y/{\cal Y}$ corresponding to $W$.}

\bigskip\noindent
{\it 2.D $T$-labelling of $\cal Y$.}

\smallskip
In this rather long subsection we introduce the concept of a $T$-labelling 
of $\cal Y$, which is the most important tool in our proof of Proposition 2.1.
Recall that
$T$ denotes the countable tree with infinite valence at every vertex. We fix terminology and notation concerning various objects inside $T$. 
We choose a {\it base vertex} in $T$, denoting it $t_0$. A {\it central ray} in $T$ is any infinite path $\gamma$ started at $t_0$, with consecutive vertices denoted 
$\gamma(0),\gamma(1),\dots$. For any vertex $t\in V_T\setminus\{ t_0 \}$ its {\it ancestor} $a_t$ is the adjacent vertex on the path from $t$ to $t_0$. A {\it sector}  based at $t$, denoted $\Sigma_t$, is the set of all $s\in V_T$ for which $t$ lies on the path from $s$ to $t_0$ (including $s=t$); $t$ is then called {\it the base} of the sector $\Sigma_t$.  The set of {\it succesors} of $t$ is the set $\Sigma_t^1=\{s\in\Sigma_t:d_T(s,t)=1\}$. For any integer $k\ge0$ the $k$-{\it ball} $B_k$ and the $k$-{\it sphere} $S_k$ are defined as $B_k=\{ t\in V_T:d_T(t,t_0)\le k \}$, $S_k=\{ t\in V_T:d_T(t,t_0)=k \}$. 

\medskip\noindent
{\bf 2.D.1 Definition.}
Given an $(X)$-regularizing family $\cal Y$ of subspaces in a metric compactum $Y$, 
a {\it $T$-labelling} for $\cal Y$ is a labelling $(Y_t)_{t\in V_T}$ of $\cal Y$ by elements of the set $V_T$ such that:
\itemitem{(L1)} the map $t\to Y_t$ is a bijection from $V_T$ to $\cal Y$; 
\itemitem{(L2)} for any central ray $\gamma$ in $T$ the sequence of subspaces $Y_{\gamma(n)}$ converges to a point in the complement $Y^\bullet$;
\itemitem{(L3)} for each $t\in V_T\setminus\{ t_0 \}$ the family $Y_{s}:s\in\Sigma_t^1$ approximates the subspace $Y_t$; similarly, the family $Y_t:t\in S_1$ approximates the subspace $Y_{t_0}$;
\itemitem{(L4)} for each $t\in V_T\setminus\{ t_0 \}$, closure in $Y$ of the union of the family $\{ Y_s:s\in\Sigma_t \}$, denoted $H_t$, is a closed and open subset of $Y$ which is disjoint with $Y_{a_t}$; 
\itemitem{(L5)} $\lim_{t\ne t_0}\hbox{diam}(H_t)=0$;
\itemitem{(L6)} for any two distnct $t_1,t_2\in S_1$, as well as for any
$s\in V_T\setminus\{ t_0 \}$ and any two distinct $t_1,t_2\in\Sigma_s$,
we have $H_{t_1}\cap H_{t_2}=\emptyset$.

\medskip
The next result shows that a $T$-labelling is a potentially useful tool for proving Proposition 2.1.

\medskip\noindent
{\bf 2.D.2 Proposition.}
{\sl Let $X$ be a nonempty metric compactum, and let $Y$ be an $(X)$-regular space, with $(X)$-regularizing family $\cal Y$ of subspaces. If $\cal Y$ admits a $T$-labelling then $Y$ is homeomorphic to the dense amalgam $\widetilde\sqcup( X)$.}

\medskip\noindent
{\bf Proof:}
Let $(Y_t)_{t\in V_T}$ be a $T$-labelling for $\cal Y$.

\medskip\noindent
{\it Step 1. A complete tree system compatible with the $T$-labelling.}

We start with constructing a complete tree system of peripheral extensions
for $X$, $\Theta=(\{ X_t:t\in V_T \},\{ b_t:t\in V_T \})$, satisfying the following conditions:

\item{(1)} for each $t\in V_T$ we have $X_t=Y_t$;
\item{(2)} for each $t\in V_T$ the map $b_t:N_t\to P_t=\ddot{X}_t\setminus X_t$ satisfies the following: choosing any points $\beta_t(s)\in Y_t$ such that 
$d_Y(\beta_t(s),Y_{s})=d_Y(Y_t,Y_{s})$, for all $s\in N_t, $we have
$$
\lim_{s\in N_t}d_{\ddot{X}_t}(b_t(s),\beta_t(s))=0 
. \leqno{(2.{\rm D}.2.1)}
$$

To construct maps $b_t$ satisfying (2.D.2.1), we proceed for each $t\in V_T$ indpendantly as follows. Order the vertices of $N_t$ into a sequence $(s_n)$ and the points of $P_t$ into a sequence $(x_n)$. Recall that, by condition (L3) in Definition 2.D.1, we have $\lim_{n}d_Y(\beta_t(s_n),Y_{s_n})=0$, and thus for any $y\in Y_t$
there is a subsequence $n_m$ such that $\lim_m\beta_t(s_{n_m})=y$. Iterate the following two steps, starting with $n=1,2$. For odd $n$, if $j$ is the smallest index for which $\beta_t(s_{j})$ has not yet been defined,
put $b_t(s_{j})=x$ for any $x\in P_t$ which was not yet chosen as the image of any other $s$, and which satisfies $d_{\ddot{X}_t}(x,\beta_t(s_{j}))<{1\over n}$. For even $n$, if $j$ is the smallest idex for which $x_j$ has not yet been chosen as the image of any $s$, choose any $s$ for which $b_t(s)$ has not yet been defined and such that $d_{\ddot{X}_t}(x_j,\beta_t(s))<d_{\ddot{X}_t}(x_j,Y_t)+{1\over n}$. We skip the direct verification that $b_t$ is then a bijection and satisfies (2.D.2.1).

\medskip\noindent
{\it Step 2: the map $h:\vee\Theta\setminus P_\Theta\to Y$.}

Recall that we have the identification 
$\vee\Theta\setminus P_\Theta=(\bigsqcup_{t\in V_T} Y_t)\sqcup\partial T$. 
If $x\in\partial T$, let $\gamma_x$ be the unique central ray in $T$ 
representing $x$. Accordingly with the above identification, 
put:
\item{$\bullet$} $h(x):=x$ if $x\in Y_t$ for some $t\in V_T$;
\item{$\bullet$} $h(x):=\lim_n Y_{\gamma_x(n)}$ if $x\in\partial T$.

\noindent
Note that, due to condition (L2), the latter limit exists and
is a point of $Y^\bullet$.

In the next three steps we will show that $h$ is respectively injective, 
surjective and open, thus getting that it is a homeomorphism. 
Since $\vee\Theta\setminus P_\Theta\cong\widetilde\sqcup(X)$, 
this will complete the proof of Proposition 2.D.2.

\medskip\noindent
{\it Step 3: $h$ is injective.} 

Since $h$ maps the subset $\bigsqcup_{t\in V_T} Y_t\subset\vee\Theta\setminus P_\Theta$ 
injectively on the subset $\cup_{t\in V_T}Y_t\subset Y$, 
and since by condition (L2) it maps $\partial T$ to the subset 
$Y^\bullet=Y\setminus(\cup_{t\in V_T}Y_t)$, 
it is sufficient to show that the restriction of $h$ to $\partial T$ is injective.

Consider two distinct points $p,q\in\partial T$, 
and the corresponding central rays $\gamma_p,\gamma_q$. 
Let $i$ be the smallest number such that $\gamma_p(i)\ne\gamma_q(i)$. 
Denote by $H_p,H_q$ respectively the closures in $Y$ of the unions 
$\cup\{ Y_s:s\in\Sigma_{\gamma_p(i)} \}$, $\cup\{ Y_s:s\in\Sigma_{\gamma_q(i)} \}$. 
By condition (L6)
we get $H_p\cap H_q=\emptyset.$ The assertion then follows by observing that 
$h(p)\in H_p$ and $h(q)\in H_q$.

\medskip\noindent
{\it Step 4: $h$ is surjective.}

Obviously, any point $x\in\cup_{t\in V_T}Y_t\subset Y$ is in the image of $h$. 
Thus, we need to show that any point $q\in Y^\bullet$ is also in this image.

According to Fact 2.B.1, there is a sequence ${t_n}$ such that in $Y$ 
we have $q=\lim_n Y_{t_n}$. Recall that for each $t\in V_T\setminus\{ t_0 \}$ 
we denote by $H_t$ the closure in $Y$ of the union $\cup\{ Y_s:s\in\Sigma_t \}$.
We claim that there is $u\in S_1$ such that 
$t_n\in\Sigma_u$ for infinitely many $n$.
Indeed, if there is no such $u$ then, denoting by $u_n$ this vertex of $S_1$ 
for which $t_n\in\Sigma_{u_n}$, we have $\lim_n\hbox{diam}(H_{u_n})=0$ 
(due to condition (L5)), and since $Y_{t_n}\subset H_{u_n}$, 
it follows that $q=\lim_n H_{u_n}$. Consequently, we also have $q=\lim Y_{u_n}$, 
and due to condition (L3) this imlies that $q\in Y_{t_0}$, 
despite $q\in Y^\bullet$. Thus, there is $u\in S_1$ such that 
$t_n\in\Sigma_u$ for infinitely many $n$. Moreover, since then $q\in H_u$, 
and since by (L6)the subsets $H_s:s\in S_1$ are pairwise disjoint, 
it follows that $u$ as above is unique. We denote it $u_1$.

Iterating the above argument, for each natural $k$ 
we get a unique $u_k\in S_k$ such that $t_n\in\Sigma_{u_k}$ for infinitely many $n$. 
By uniqueness of $u_k$, we get that $\Sigma_{u_{k+1}}\subset\Sigma_{u_k}$ 
for each $k$, and thus the sequence $t_0,u_1,u_2,\dots$is a central ray in $T$. 
Denote by $p\in\partial T$ the point corresponding to this central ray. 
Since we have $q\in H_{u_k}$ for each $k$, we also have $q=\lim_k H_{u_k}$, 
and consequently $q\in\lim_k Y_{u_k}=h(p)$. This completes the proof 
of surjectivity.


\medskip\noindent
{\it Step 5: $h$ is open.}

We refer to the basis ${\cal D}$ of the topology of 
$\vee\Theta\setminus P_\Theta=(\bigsqcup_{t\in V_T} Y_t)\sqcup\partial T$, 
as described in Proposition 1.8. 
We need to show that far any set $D(t,U)\in{\cal D}$
(where $t\in V_T$ is a vertex, and $U\subset\ddot{Y}_t$ is an open subset)
its image $h(D(t,U))$ 
is an open subset of $Y$.

Recall that there are three kinds of points in $D(t,U)$:
\item{(1)} points $x\in U\cap Y_t$;
\item{(2)} points $y\in Y_s$ for $s\ne t$ such that $b_t([t,t_s])\in U$;
\item{(3)} points $p\in\partial T$ such that $b_t([t,t_z])\in U$.

\noindent
We will show that the image $z$ of a point of each kind is contained in $h(G_U)$ 
together with some open neighbourhood of $z$ in $Y$.

Let $z=h(p)$ for some $p$ of kind (3) above. 
Choose any vertex $s\in V_T\setminus\{ t_0 \}$ lying on the central ray 
from $t_0$ to $p$ and such that $t\notin\Sigma_s$. 
Note that then the set $\sqcup\{ Y_u:u\in\Sigma_s \}\i D(t,U)$ 
and the set of all $q\in\partial T$ represented by central rays 
passing through $s$ are both the subsets of $D(t,U)$. 
We also claim that, denoting the union of these two subsets by $D_s$,
we have $h(G_s)=H_s$. The inclusion $h(G_s)\i H_s$ is obvious. 
For the opposite inclusion, the argument is the same as that in Step 4. 
Thus, we get $z\in H_s=\i h(D(t,U))$, 
where the subset $H_s$ is open (by condition (L4)).

Now, let $z=h(y)$ for some $y\in Y_s$ of kind (2) above. 
We consider three subcases concerning the position of $s$. 
First, suppose that $s$ is not lying on the path from $t_0$ to $t$. 
Then, arguing as in the previous case, we get similarly that $z\in H_s\i h(D(t,U))$. 
In the remaining cases, denote by $s'$ the vertex adjacent to $s$ 
on the path from $s$ to $t$. If $s=t_0$, one shows similarly 
(using the fact that $h$ is a bijection) that 
$z\in Y\setminus H_{s'}\i h(D(t,U))$. Since by (L4) the set $H_{s'}$ is closed, 
its complement $Y\setminus H_{s'}$ is open, 
and thus it is as required. Finally, if $s$ lies in the interior 
of the path from $t_0$ to $t$, 
by condition (L4) we have $Y_s\cap H_{s'}=\emptyset$.
We then get $z\in H_s\setminus H_{s'}\i h(D(t,S))$, 
where $H_s\setminus H_{s'}$ is easily seen to be open, again due to (L4).

In the last case, let $z=h(x)$ for some $x\in U\cap Y_t$ 
(i.e. $x$ is of kind (1) above). 
Since $U$ is open, there is $\epsilon>0$ such that 
$$d_{\ddot{Y}_t}(x,b_t(s))>\epsilon \hbox{\quad for each \quad} 
s\in N_t\setminus b_t^{-1}(U). \leqno{(2.{\rm D}.2.2)}$$
In view of (2.D.2.1), we then have 
$$\liminf_{s\in N_t\setminus 
b_t^{-1}(U)}d_{\ddot{Y}_t}(x,\beta_t(s))\ge\epsilon.$$
Since the metrics $d_{\ddot{Y}_t}$ and $d_Y$ restricted to 
$Y_t$ are equivalent, and since 
$x\in Y_t\subset(\bigsqcup_{t\in V_T}Y_t)\sqcup\partial T$ 
coincides with $z=h(x)\in Y_t\i Y$, there is $\epsilon'>0$ such that
$$ \liminf_{s\in N_t\setminus b_t^{-1}(U)}d_Y(z,\beta_t(s))\ge\epsilon'. $$
Since $\lim_{s\in N_t}d_Y(\beta_t(s),Y_{s})=0$ and 
$\lim_{s\in N_t}\hbox{diam}(Y_{s})=0$, it follows that
$$  \liminf_{s\in N_t\setminus b_t^{-1}(U)}d_Y(z,Y_{s})\ge\epsilon'. 
\leqno{(2.{\rm D}.2.3)}$$

For each $s\in N_t$ consider the half-tree $\Psi_s$ in $T$ containing $s$ 
and not containing $t$. Put $\Omega_s=\overline{\cup\{ Y_u:u\in\Psi_s \}}$, 
where the closure is taken in $Y$. Note that that for all $s\in N_t$ 
except possibly one (namely this for which $t_0\in\Psi_s$) we have 
$\Omega_s=H_s$, and hence 
$$
\lim_{s\in N_t}\hbox{diam}(\Omega_s)=0. 
\leqno{(2.{\rm D}.2.4)}
$$
Since $Y_s\i \Omega_s$ for each $s\in N_t$, 
it follows from (2.D.2.3) and (2.D.2.4) that
$$
\liminf_{s\in N_t\setminus b_t^{-1}(U)}d_Y(z,\Omega_s)\ge\epsilon'. 
$$
Thus, 
$$
\hbox{for almost all }s\in N_t\setminus b_t^{-1}(U) 
\hbox{ we have } d_Y(z,\Omega_s)>{\epsilon'\over2}. 
\leqno{(2.{\rm D}.2.5)} 
$$
We claim also that for any $s\in N_t$ we have $d_Y(z,\Omega_s)>0$. To see this, 
it is enough to note that for each $s\in N_t$ we have $Y_t\cap \Omega_s=\emptyset$.
Indeed, this is true by condition (L4) for all $s$ except possibly this one 
for which $t_0\in\Psi_s$. We denote this exceptional $s$ by $s_0$. 
If this $s_0$ exists,  
one easily notes that, since by (L4) the subset $H_t$ is open, 
we have $H_t\cap \Omega_{s_0}=\emptyset$, and consequently 
$Y_t\cap\Omega_{s_0}=\emptyset$. 

As a consequence of the assertions in the previous paragraph, 
there is $\delta>0$ such that $d_Y(z,\Omega_s)>\delta$ for all 
$s\in N_t\setminus b_t^{-1}(U)$ and $d_Y(z,Y_t\setminus U)>\delta$. 
Since from the definition of $h$ one deduces easily that
$$
h([\vee\Theta\setminus P_\Theta]\setminus D(t,U)
\i(Y_t\setminus U)\cup\bigcup\{ \Omega_s:s\in N_t\setminus b_t^{-1}(U) \},
$$ 
it follows from bijectivity of $h$ that the metric ball $B_\delta(x,(Y,d_Y))$ is contained in $h(D(t,U))$.
This completes the proof of openness of $h$, 
and hence also the proof of Proposition 2.D.2.



\bigskip\noindent
{\bf Proof of Proposition 2.1.}

\smallskip
Let $Y$ be an $(X)$-regular space, with $(X)$-regularizing family $\cal Y$.
In view of Proposition 2.D.2, to prove Proposition 2.1,  
it is suficient to show that $\cal Y$ admits a $T$-labelling.
Before starting the actual construction of such a $T$-labelling,  
order $\cal Y$ into a sequence $(Y_k)_{k\ge1}$.
We demand that a labelling that we construct satisfies the following:

\itemitem{(r1)} for each $k\ge1$ we have $Y_k\in\{ Y_u:u\in B_{k-1} \}$;
\itemitem{(r2)} for each $k\ge1$ and any $t\in S_k$, denoting by $[t_0,u_1,\dots,u_{k-1},t]$ the path in $T$ from $t_0$ to $t$, and putting
$$
d_t:=\min\lgroup \hbox{diam}(Y_t),{1\over3}d_Y(Y_k,Y_t),{1\over3^2}d_Y(Y_{k-1},Y_{u_{k-1}}),\dots,{1\over3^k}d_Y(Y_1,Y_{u_{1}})\rgroup,
$$
for any $s\in\Sigma_t^1$ we have $Y_s\i N_{d_t}(Y_t)$ and $\hbox{diam}(Y_s)<{1\over2}\hbox{diam}(Y_t)$.

\medskip
Note that, in view of Lemma 2.B.3, we have the following.

\medskip\noindent
{\bf Claim.}
{\sl If a labelling $(Y_t)_{t\in V_T}$ for $\cal Y$ satisfies the above conditions (r1) and (r2) then it satisfies condition (L2) of Definition 2.D.1.}

\medskip
We start the inductive construction of a $T$-labelling for $\cal Y$ by putting $Y_{t_0}:=Y_1$. Induction proceeds with respect to radii of
balls $B_n$ and spheres $S_n$ in $V_T$. At the first essential
(i.e. not trivial) step, for each $t\in S_1$ we choose $Y_t$ so that
\itemitem{(1)} the family $\{ Y_t:t\in S_1 \}$ contains $Y_2$ and approximates $Y_{t_0}$;
\itemitem{(2)} if we put $d_t=\min\lgroup\hbox{diam}(Y_t),{1\over3}d_Y(Y_1,Y_t)\rgroup$, then the family 
$N_{d_t}(Y_t):t\in S_1$ covers $Y\setminus Y_{t_0}$;
\itemitem{(3)}
for any $Z\in{\cal Y}\setminus\{ Y_u:u\in B_1 \}$ there is $t\in S_1$ such that $Z\i N_{d_t}(Y_t)$ and $\hbox{diam}(Z)<{1\over2}\hbox{diam}(Y_t)$.

\noindent
To make such a choice, consider the subset $E_0=\{ x\in Y:d_Y(x,Y_{t_0}) \}\ge1$, and for each $m\ge1$ consider the subset $E_m=\{ x\in Y:2^{-m}\le d_Y(x,Y_{t_0})\le2^{-m+1} \}$. Each of those subsets is closed in $Y$, and hence compact. For each $W\in{\cal Y}\setminus\{  Y_{t_0} \}$ put $d_W=\min\lgroup\hbox{diam}(W),{1\over3}d_Y(Y_1,W)\rgroup$. For each $m\ge0$  choose a finite subfamily ${\cal W}_m\i{\cal Y}\setminus\{ Y_{t_0} \}$ such that each $W\in{\cal W}_m$ intersects $E_m$, and
the corresponding family of neighbourhoods $\{ N_{d_W}(W):W\in{\cal W}_m \}$ covers $E_m$. 
Denote by ${\cal W}_m^+$ the set of all $W'\in{\cal Y}$, $W'\ne Y_{t_0}$, $W'\cap E_m\ne\emptyset$, such that 
$W'$ is not contained in any single neighbourhood from the family $\{ N_{d_W}(W):W\in{\cal W}_m \}$ or $\hbox{diam}(W')\ge{1\over2}\min\lgroup \hbox{diam}(W):W\in{\cal W}_m \rgroup$. Note that for each $m\ge0$ the family ${\cal W}_m^+$ is finite. Put ${\cal W}:=[\bigcup_{m\ge0}({\cal W}_m\cup{\cal W}_m^+)]\cup\{ Y_2 \}$ and label $\cal W$ using $S_1$ as the set of labels, so that ${\cal W}=\{ Y_t:t\in S_1 \}$. Observe that conditions (1)--(3) above are then satisfied
(we skip a rather straightforward argument).

By Proposition 2.A.1, the space $(Y/{\cal Y})\setminus\{[Y_{t_0}]\}$ is homeomorphic to the punctured Cantor space $C_0$
(here, for $s\in V_T$ we denote by $[Y_{s}]$ the point of $Y/{\cal Y}$
corresponding to $Y_{s}$). 
Moreover, since the just chosen family $Y_t:t\in S_1$ approximates $Y_{t_0}$, the corresponding subset $\{ [Y_t]:t\in S_1 \}$ is discrete in $(Y/{\cal Y})\setminus\{[Y_{t_0}]\}$ (see Fact 2.C.2(b)). 
For each $t\in S_1$ put
$$
U_t:=N_{d_t}(Y_t)\setminus 
\bigcup\{ Z\in{\cal Y}: Z\not\subset N_{d_t}(Y_t) \}
\setminus 
\bigcup\{ Z\in{\cal Y}:Z\ne Y_t,\hbox{diam}(Z)\ge
{1\over2}\hbox{diam}(Y_t) \}.
$$
Observe that nullness of $\cal Y$ has the following consequences.
First, the union $\bigcup\{ Z\in{\cal Y}:Z\ne Y_t,\hbox{diam}(Z)\ge
{1\over2}\hbox{diam}(Y_t) \}$ is finite, and hence it yields a closed subset of $Y$. 
Second, the set
$$
\overline{\bigcup\{ Z\in{\cal Y}: Z\not\subset N_{d_t}(Y_t) \}}
\setminus 
\bigcup\{ Z\in{\cal Y}: Z\not\subset N_{d_t}(Y_t) \}
$$
(where the closure is taken in $Y$) is disjoint with $N_{d_t}(Y_t)$.
It follows that
$$
U_t=N_{d_t}(Y_t)\setminus 
\overline{\bigcup\{ Z\in{\cal Y}: Z\not\subset N_{d_t}(Y_t) \}}
\setminus 
\bigcup\{ Z\in{\cal Y}:Z\ne Y_t,\hbox{diam}(Z)\ge
{1\over2}\hbox{diam}(Y_t) \}.
$$
In particular, $U_t$ is an open neighbourhood of $Y_t$ in $Y$.
Moreover, by conditions (2) and (3) above, the family $U_t:t\in S_1$
is a covering of $Y\setminus Y_{t_0}$. Obviously, the sets $U_t$
are all $\cal Y$-saturated. Thus, their images $U_t'$ through the
quotient map $Y\to Y/{\cal Y}$ form an open covering of 
$(Y/{\cal Y})\setminus\{[Y_{t_0}]\}$ by the sets whose
closures in $(Y/{\cal Y})\setminus\{[Y_{t_0}]\}$ are compact
(because their closures in $Y/{\cal Y}$ do not contain the point
$[Y_{t_0}]$), and for each $t\in S_1$ we have $[Y_t]\in U_t'$.
By Proposition 2.A.4, there is a partition of the space 
$(Y/{\cal Y})\setminus\{[Y_{t_0}]\}$ into subsets $K_t:t\in S_1$ 
which are compact, open, and
such that for each $t$ we have $[Y_t]\in K_t\i U_t'$. Denoting by $q:Y\setminus Y_{t_0}\to (Y/{\cal Y})\setminus\{[Y_{t_0}]\}$ the quotient map, we get the partition of $Y\setminus Y_{t_0}$ into subsets $L_t=q^{-1}(K_t):t\in S_1$ which are closed and open in $Y$ and $\cal Y$-saturated. It is not hard to see that for each $t\in S_1$ we also have 
\itemitem{(p1)}
$Y_t\i L_t$;
\itemitem{(p2)} $L_t\i N_{d_t}(Y_t)$ and thus, since $d_t\le\hbox{diam}(Y_t)$, we have $\hbox{diam}(L_t)<3\hbox{diam}(Y_t)$;
\itemitem{(p3)} each $Z\in {\cal Y}$ contained in $L_t$ 
and distinct from $Y_t$ satisfies 
$\hbox{diam}(Z)<{1\over2}\hbox{diam}(Y_t)$.

\smallskip
We now proceed to the general inductive step of the construction.
Suppose that for some $n\ge1$ and for all $t\in B_n$ we have already 
chosen the subspaces $Y_t$ so that the family $Y_t:t\in B_n$
contains all of the subspaces $Y_1,Y_2,\dots,Y_{n+1}$.
Suppose also that we have constructed a partition of the subspace
$Y\setminus\cup\{ Y_u:u\in B_{n-1} \}$ into $\cal Y$-saturated
subspaces $L_t:t\in S_n$, each open and closed in $Y$, such that
for each $t\in S_n$ we have
\itemitem{(p1*)}
$Y_t\i L_t$;
\itemitem{(p2*)} $L_t\i N_{d_t}(Y_t)$ and thus $\hbox{diam}(L_t)<3\hbox{diam}(Y_t)$;
\itemitem{(p3*)} each $Z\in {\cal Y}$ contained in $L_t$ 
and distinct from $Y_t$ satisfies 
$\hbox{diam}(Z)<{1\over2}\hbox{diam}(Y_t)$.

\noindent
For each $t\in S_1$ do the following. For each $s\in\Sigma_t^1$ choose $Y_s$ so that
\itemitem{(t1)} $Y_s\i L_t\setminus Y_t$ (then $Y_s\i N_{d_t}(Y_t)$ and $\hbox{diam}(Y_s)<{1\over2}\hbox{diam}(Y_t))$;
\itemitem{(t2)} if $Y_{n+2}\i L_t\setminus Y_t$ then for some $s\in\Sigma_t^1$ we have $Y_s=Y_{n+2}$;
\itemitem{(t3)} the family $\{ Y_s:s\in\Sigma_t^1 \}$ approximates $Y_t$;
\itemitem{(t4)} denoting by $[t_0,u_1,\dots,u_{n},s]$
the path in $T$ from $t_0$ to $s$,
if we put 
$$
d_s=\min\lgroup \hbox{diam}(Y_s), {1\over3}d_Y(Y_{n+1},Y_s), {1\over3^2}d_Y(Y_n,Y_{u_n}), 
\dots,{1\over3^{n+1}}d_Y(Y_1,Y_{u_1}) \rgroup,
$$ 
then 
the family $N_{d_s}:s\in\Sigma^1_t$ covers $L_t\setminus Y_t$ and for any $Z\in{\cal Y}$, $Z\i L_t\setminus Y_t$, there is $s\in\Sigma_t^1$ such that $Z\i N_{d_s}(Y_s)$ and 
$\hbox{diam}(Z)<{1\over2}\hbox{diam}(Y_s)$.

\noindent
This can be done in a way similar as described above (for the family $Y_t:t\in S_1$). Moreover, since the quotient space $L_t/\{ Z\in{\cal Y}:Z\i L_t \}$ 
is homeomorphic to the Cantor space (because it is an open and closed subspace of $Y/{\cal Y}$), arguing as above with the help of 
Proposition A.2.4, we obtain a partition of
$L_t\setminus Y_t$ into a family of subsets $L_s:s\in\Sigma_t^1$ which are closed and open in $L_t$ (and hence also in $Y$), $\cal Y$-saturated, and such that $Y_s\i L_s\i N_{d_s}(Y_s)$ for each $s\in\Sigma^1_t$. 
Moreover, for each such $s$ condition (p3*) holds with $s$ substituted for $t$.

By the above  inductive construction, we get an injective map $u\to Y_u$ from $V_T$ to $\cal Y$, which is also surjective due to condition (r1) (which follows from condition (t2)). This map is thus a labelling of $\cal Y$, i.e. condition (L1) of Definition 2.D.1 holds. Moreover, 
since conditions (r1) and (r2) are fulfilled due to (t2) and (t4),
Claim above ensures that this labelling satisfies condition (L2).  
Condition (L3) holds by (t3).
It follows fairly directly from (t1) and from the construction and properties of sets $L_s$ that for each $t\ne t_0$ we have $H_t\i L_t$.
This easily implies condition (L6), and in view of (p2*) it also
implies (L5).
The argument as in Step 4 of the proof of Proposition 2.D.2 shows that
in fact  for each $t\ne t_0$ the subspace $H_t$ coincides with $L_t$, thus being open and closed in $Y$, which justifies condition (L4). Hence, the construction above yields a $T$-labelling for $\cal Y$,
which completes the proof.

\bigskip
\noindent
{\bf 3. Properties of the operation of dense amalgam.}

\medskip
In this section we derive, using the construction and the characterization,
various properties of the operation of dense amalgam. In particular,
we justify Proposition 0.1 of the introduction.

\smallskip
Start with observing that part (1) of Proposition 0.1 (which yields
commutativity of the dense amalgam) follows by the definition
of $\widetilde\sqcup(X_1,\dots,X_k)$,
as given at the end of Section 1.

We next pass to part (3) of Proposition 0.1. We will show the following
result, which obviously implies the statement of part (3),
but in fact it is easily seen to be just equivalent.

\medskip\noindent
{\bf 3.1 Proposition.}
{\it Let $X_1,\dots,X_k$ be any nonempty metric compacta, and let
$M_1,\dots,M_m$ be the 
pairwise non-homeomorphic spaces
representing all homeomorphism types appearing among $X_1,\dots,X_k$.
Then $\widetilde\sqcup(X_1,\dots,X_k)=\widetilde\sqcup(M_1,\dots,M_m)$.}

\medskip\noindent
{\bf Proof:}
In view of the characterization given in Theorem 0.2, and proved in
Section 2, it is sufficient to show that the space 
$Y=\widetilde\sqcup(X_1,\dots,X_k)$ is $(M_1,\dots,M_m)$-regular.
Let $\cal Y$ be an $(X_1,\dots,X_k)$-regularizing family for $Y$.
It is partitioned into subfamilies ${\cal Y}_i:1\le i\le k$ so that conditions
(a1)--(a5) of the introduction hold. For each $i\in\{ 1,\dots,k \}$ let $m(i)$
be this index for which $X_i$ is homeomorphic to $M_{m(i)}$.
We define a new partition of $\cal Y$, into subfamilies
${\cal Y}_j':1\le j\le m$, putting ${\cal Y}_j':=\cup\{ {\cal Y}_i:m(i)=j \}$
for each $j\in\{ 1,\dots,m \}$. A direct verification shows that
$\cal Y$ equipped with this new partition is an 
$(M_1,\dots,M_m)$-regularizing family for $Y$, which completes the proof.

\medskip
To deal with the remaining parts of Proposition 0.1, we will need
the following auxilliary result, which follows fairly directly from
the construction of the dense amalagam, as described in Section 1.

\medskip\noindent
{\bf 3.2 Lemma.}
{\it Suppose that $Y$ is an $(X_1,\dots,X_k)$-regular space, with
$(X_1,\dots,X_k)$-regulari\-zing family $\cal Y$, and let $Z\in{\cal Y}$.
Suppose that $Z$ is not connected, and let $A_1,A_2$ be nonempty 
open and closed subsets of $Z$ forming its partition. 
Then there is a partition
of $Y$ into open and closed subsets $H_1,H_2$ such that:}
\item{(1)} {\it $A_1=H_1\cap Z$ and $A_2=H_2\cap Z$;}
\item{(2)} {\it each subset $W\in{\cal Y}\setminus\{ Z \}$ is contained
either in $H_1$ or in $H_2$.}

\medskip\noindent
{\bf Proof:}
We may assume that $Y=\widetilde\sqcup(X_1,\dots,X_k)$, and we identify
it with $(\bigsqcup_{t\in V_T}X_t)\sqcup\partial T$, as in Section 1,
where each $X_t$ is homeomorphic to $X_1\sqcup\dots\sqcup X_k$.
Under this homeomorphism, we express each $X_t$ as 
$X_t=X_{1,t}\sqcup\dots\sqcup X_{k,t}$, where $X_{i,t}\cong X_i$.
We then identify $Z$ with a subset $X_{i,t}$, for some $t\in V_T$
and some $1\le i\le k$.
We extend the partition of $Z=X_{i,t}$ into subsets $A_1,A_2$ 
first to a partition of $X_t$ into $A_1'=A_1$ and $A_2'=X_t\setminus A_1$, and then
to a partition of the peripheral extension $\ddot{X}_t$ into
open and closed subsets $A_1'',A_2''$. Putting 
$H_i:=D(t,A_i'')$, for $i=1,2$, we get a partition of $Y$ as required,
which completes the proof.

\medskip
We now turn to proving parts (2), (4) and (5) of  Proposition 0.1.

\medskip\noindent
{\bf Proof of Proposition 0.1(2).}

In view of the characterization of dense amalgams given in Theorem 0.2,
it is sufficient to show that the space 
$Y=\widetilde\sqcup(X_1,\dots,X_{i-1},\widetilde\sqcup(X_i,\dots, X_k))$
is $(X_1,\dots,X_k)$-regular. We start with describing a natural
candidate for an $(X_1,\dots,X_k)$-regularizing family ${\cal Y}^*$.
We refer to the identification (as presented in Section 1) 
of $Y$ with the space
$(\bigsqcup_{t\in V_T}X_t)\sqcup\partial T$,
where each $X_t$ is homeomorphic to the space
$X_1\sqcup\dots\sqcup X_{i-1}\sqcup\Omega$,
with $\Omega=\widetilde\sqcup(X_i,\dots,X_k)$.
We realize each such homeomorphism as 
$X_t=X_1^t\sqcup\dots\sqcup X_{i-1}^t\sqcup\Omega^t$.
We also identify each $\Omega^t$ with the space 
$(\bigsqcup_{s\in V_{T_t}}X_{t,s})\sqcup\partial T_t$,
where $T_t$ is a separate copy of the tree $T$, and 
where each $X_{t,s}$ is homeomorphic to $X_i\sqcup\dots\sqcup X_k$,
which we write as $X_{t,s}=X_i^{t,s}\sqcup\dots\sqcup X_k^{t,s}$.
Now, for each $j\in\{ 1,\dots,i-1 \}$ we put
${\cal Y}_j^*:=\{ X_j^t:t\in V_T \}$, 
and  for each $j\in\{ i,i+1,\dots,k \}$ we put
${\cal Y}_j^*:=\{ X_j^{t,s}:t\in V_T, s\in V_{T_t} \}$.

Verification that the so described family 
${\cal Y}^*=\bigsqcup_{j=1}^k{\cal Y}_j^*$ satisfies conditions (a1)--(a4)
of an $(X_1,\dots,X_k)$-regularizing family is straightforward, and we 
skip it. Verification of condition (a5) is a bit more involved.
Let $x,y$ be two points of $Y$ not belonging to the same set of
${\cal Y}^*$. We need to separate $x$ from $y$ by a 
${\cal Y}^*$-saturated open and closed subset of $Y$.
To do this, one needs to consider several cases of positions of $x$ and $y$.
We verify the statement in the case when for some $t_0\in V_T$, some 
$s_0\in V_{T_t}$, and some $j\in\{ i,i+1,\dots,k \}$ we have 
$x\in X_j^{t_0,s_0}$ and $y\in\partial T_{t_0}$. (The arguments in the 
remaining cases are similar, and we omit them.)

Consider the family ${\cal Y}^{t_0}:=\{ X_{t_0,s}:s\in V_{T_{t_0}} \}$,
which is the canonical $(X_i,\dots,X_k)$-regularizing family for 
$\Omega^{t_0}$, as indicated at the end of Section 1. By condition (a5), there is a partition of $\Omega^{t_0}$ into open, closed and ${\cal Y}^{t_0}$-saturated
subsets $A_1,A_2$ such that $x\in A_1$ and $y\in A_2$.
Consider then the family ${\cal Y}:=
(\bigcup_{j=1}^{i-1}\{ X_j^t:t\in V_T \})\cup
\{ \Omega^t:t\in V_T \}$,
which is the canonical $(X_1,\dots,X_{i-1},\Omega)$-regularizing
family for $Y$. By Lemma 3.2, there is a partition of $Y$ into open
and closed subsets $H_1,H_2$ such that $A_i=H_i\cap\Omega^{t_0}$
for $i=1,2,$ and each subset of ${\cal Y}\setminus\{ \Omega^{t_0} \}$
is contained either in $H_1$ or in $H_2$.
Since the subsets $H_i$ are clearly ${\cal Y}^*$-saturated,
and since $x\in H_1$ and $y\in H_2$, the assertion follows in the considered case, which completes the proof.

\vfill\break

\medskip\noindent
{\bf Proof of Proposition 0.1(4).}

Let $Y=\widetilde\sqcup(X_1,\dots,X_k,Q)$.
In view of Theorem 0.2, it is sufficient to show that $Y$ is 
$(X_1,\dots,X_k)$-regular. Let $\cal Y$ be the canonical
$(X_1,\dots,X_k,Q)$-regularizing family for $Y$, as described at the end
of Section 1. More precisely, viewing $Y$ as 
$(\bigsqcup_{t\in V_T}X_t)\sqcup\partial T$, with each $X_t$
homeomorphic to $X_1\sqcup\dots\sqcup X_k\sqcup Q$, which we write as
$X_t=X_{1,t}\sqcup\dots\sqcup X_{k,t}\sqcup Q_t$, we have
${\cal Y}=(\bigsqcup_{i=1}^k \{ X_{i,t}:t\in V_T \})
\sqcup \{ Q_t:t\in V_T \}$.
As a natural candidate for an $(X_1,\dots,X_k)$-regularizing
family for $Y$ we take ${\cal Y}'=\bigsqcup_{i=1}^k{\cal Y}_i'$
with ${\cal Y}_i'=\{ X_{i,t}:t\in V_T \}$.

We need to check conditions (a1)--(a5) for ${\cal Y}'$.
The only one for which the verificatioin is not starightforward is
condition (a5). The only problematic case is when we consider points
$x,y\in Q_{t_0}$ for some $t_0\in V_T$ (which no longer belong to the same set of ${\cal Y}'$). Since $Q_{t_0}$ (being homeomorphic to $Q$)
is totally disconnected, we may choose its partition into open and closed
subsets $A_1,A_2$ such that $x\in A_1$, $y\in A_2$. Then,
applying Lemma 3.2 to the family $\cal Y$ and to $Z=Q_{t_0}$,
we get partition of $Y$ into open and closed subsets $H_1,H_2$
which are ${\cal Y}'$-saturated. Since obviously we have
$x\in H_1$ and $y\in H_2$, the proof is completed.

\medskip\noindent
{\bf Proof of Proposition 0.1(5).}

We refer to the characterization of the Cantor space $C$ as the compact metric space which is totally disconnected and has no isolated points.
Since the arguments are standard and similar to the previous ones, we only sketch them.

The malgam $\widetilde\sqcup(Q)$ is compact and metrizable by the
argument provided in Section 1 for all dense amalgams.
It has no isolated points by conditions (a3) and (a4) (this is again true
for any dense amalgam). Finally, $\widetilde\sqcup(Q)$ is totally
disconnected due to condition (a5), and by total disconnctedness of $Q$
combined with Lemma 3.2. We omit further details.

\bigskip

\noindent
{\bf 4. $E{\cal Z}$-boundaries for graphs of groups.}

\medskip
This section is devoted to the proof of Theorem 0.3(1).
More precisely, given a non-elementary graph of groups $\cal G$ with finite edge groups,
and with vertex groups equipped with $E{\cal Z}$-boundaries 
$\partial G_v$,
we show that the model of $E{\cal Z}$-boundary for the fundamental
group of $\cal G$ constructed by Alexandre Martin in [Ma]
is homeomorphic to the dense amalgam of the boundaries 
$\partial G_v$.

\bigskip\noindent
{\bf 4.1. Graphs of groups.}

\smallskip
We recall basic terminology and notation concerning graphs of groups,
referring the reader to [Se] for a more complete exposition.
We consider graphs $Y$ with multiple edges and loop edges allowed.
We denote by $V_Y$ the set of vertices, and by $O_Y$ the set of oriented edges
of $Y$. Given $a\in O_Y$, we denote by $\alpha(a)$ and $\omega(a)$
the initial and the terminal vertex of $a$, respectively.
For $a\in O_Y$, we denote by $\bar a$ the oppositely oriented edge,
and by $|a|$ the nonoriented edge underlying $a$. The set of nonoriented
edges of $Y$ will be denoted $|O|_Y$.

\medskip\noindent
{\bf 4.1.1 Definition.} A {\it graph of groups} over a graph $Y$ is a tuple
$$
{\cal G}=(\{ G_v:v\in V_Y \}, \{ G_e:e\in |O|_Y \}, \{ i_a:a\in O_Y\}),
$$
where $G_v$ and $G_e$ are groups, and $i_a:G_{|a|}\to G_{\omega(a)}$
are group monomorphisms.

\medskip
Given a graph $Y$, we denote by $Y'$ its first barycentric subdivision.
For any $a\in O_Y$, we denote by $a^+$ the nonoriented edge in $Y'$
which connects the barycenter of  $|a|$ with the vertex $\omega(a)$.
Thus, the set of nonoriented edges of $Y'$ is exactly
$\{ a^+:a\in O_Y \}$. 

\medskip\noindent
{\bf 4.1.2 Definition.}
Let $\cal G$ be a graph of groups over a graph $Y$, and let $\Xi$ be
a maximal tree in $Y'$. Consider the set of symbols 
$S=\{ s_a: a\in O_Y,a^+\not\subset \Xi \}$.
The {\it fundamental group} $G=\pi_1({\cal G},\Xi)$ is the group
$$
G=\big((*_{v\in V_Y}G_v)*(*_{e\in|O|_Y}G_e)*F_S\big)/N,
$$
where $F_S$ is the free group with the standard generating set $S$,
and where $N$ is the normal subgroup of the free product
$(*_{v\in V_Y}G_v)*(*_{e\in|O|_Y}G_e)*F_S$
generated by the elements
$g^{-1}i_a(g):a^+\i \Xi,g\in G_{|a|}$ and the elements
$g^{-1}s_a^{-1}i_a(g)s_a:a^+\not\subset \Xi,g\in G_{|a|}$.

\medskip
Since we have canonical injections of the groups $G_v$, $G_e$ and $F_s$
in $G$, we will often identify elements of these groups as elements of $G$.

\medskip\noindent
{\bf 4.1.3 Definition.}
Given a graph of groups $\cal G$ over $Y$, and a maximal subtree 
$\Xi\subset Y'$, {\it the Bass-Serre tree} $X=X({\cal G},\Xi)$
is described as follows:
\item{$\bullet$} $V_X=\bigsqcup_{v\in V_Y}(G/G_v)\times\{v\}$
and $O_X=\bigsqcup_{a\in O_Y}(G/G_{|a|})\times\{a\}$;

\item{$\bullet$} $\overline{(gG_{|a|},a)}=(gG_{|a|},\bar a)$;

\item{$\bullet$} $\omega((gG_{|a|},a))=
\cases{(gG_{\omega(a)},\omega(a)) & if $a^+\subset \Xi$ \cr
(gs_a^{-1}G_{\omega(a)},\omega(a)) & if $a^+\not\subset \Xi$. \cr}$

\noindent
The Bass-Serre tree $X$ comes equipped with the $G$-action (for $G=\pi_1({\cal G},\Xi)$)
given by
$$
h\cdot(gG_v,v)=(hgG_v,v) \hbox{\quad and \quad}
h\cdot(gG_{|a|},a)=(hgG_{|a|},a).
$$

\medskip
It is well known that $X=X({\cal G},\Xi)$ is 
indeed a tree, and $G$ acts on $X$
without inversions and so that the vertex and edge stabilizers are as follows:
$$
\hbox{Stab}_G((gG_v,v))=gG_vg^{-1} \hbox{\quad and \quad}
\hbox{Stab}_G((gG_{|a|},a))=gG_{|a|}g^{-1}.
$$

\noindent
There is also a canonical nondegenerate map $\pi:X\to Y$ given by
$\pi((gG_v,v))=v$ and $\pi((gG_{|a|},a))=a$, which is $G$-invariant
(i.e. $G$-equivariant with respect to the trivial action of $G$ on $Y$).

\medskip
\noindent
{\bf 4.1.4 Remark.}
A bit more geometric description of the Bass-Serre tree $X=X({\cal G},\Xi)$
(or description of its geometric realization) goes as follows.
For each $a\in O_Y$, let $\tau_{|a|}$ be
a nonoriented edge with its two associated oriented edges $\tau_a$
and $\overline{\tau_a}=\tau_{\bar a}$, and suppose that its endpoints
$\alpha(\tau_a)$ and $\omega(\tau_a)$ are distinct. View $\tau_{|a|}$
as a topological space homeomorphic to a segment. Put
$$
X=\big( \bigsqcup_{e\in|O|_Y}(G/G_e)\times\tau_e \big) / \sim,
$$
where $\sim$ is induced by the equivalences
$(gG_{|a|},\omega(\tau_a))\sim (g'G_{|a|},\omega(\tau_a))$
for the following triples $(a,g,g')\in O_Y\times G\times G$:
\item{$\bullet$} $a^+\subset \Xi$ and $g^{-1}g'\in G_{\omega(a)}$;
\item{$\bullet$} $a^+\not\subset \Xi$ and $g^{-1}g'\in s_a^{-1}G_{\omega(a)}s_a$.

\noindent
The (geometric) edges of $X$ are then the images through the quotient map of 
the relation $\sim$ of the sets $gG_e\times\tau_e$,
and we denote them $[gG_e,\tau_e]$. Similarly, the vertices of $X$ are the equivalence classes of points $(gG_{|a|},\omega(\tau_a))$,
which we denote $[gG_{|a|},\omega(\tau_a)]$.

\medskip
We now pass to discussing a not quite standard concept of a {\it non-elementary}
graph of groups, which appears in the statement of Theorem 0.3.
An oriented edge $a\in O_Y$ in 
a graph of groups 
${\cal G}=(\{G_v\},\{G_e\},\{i_a\})$ over $Y$ 
is {\it trivial} if it is not a loop
and if $i_a:G_{|a|}\to G_{\omega(a)}$ is an isomorphism. 
Given a trivial edge $a$, we define a new graph of groups ${\cal G}'$
by contracting the edge $|a|$ in $Y$ to a point (denoted $v_{|a|}$),
thus getting a new graph $Y'$, and by putting
$G'_{v_{|a|}}:=G_{\alpha(a)}$, while leaving the groups
and maps unchanged at the remaining vertices and edges.
The resulting graph of groups ${\cal G}'$ has the same fundamental group
as $\cal G$, and we say that it is obtained from $\cal G$ by an
{\it elementary collapse}. A graph of groups with no trivial edge
is said to be {\it reduced}. Obviously, any graph of groups
(over a finite graph) can be modified into a reduced graph of groups
by a sequence of elementary collapses.

\medskip
\noindent
{\bf 4.1.5 Definition.}
A graph of groups $\cal G$ over $Y$ is {\it simply elementary} if it
has one of the following three forms:
\item{$\bullet$} $Y$ consists of a single vertex, and has no edge;
\item{$\bullet$} $Y$ consists of a single vertex, $v$, and a single
loop edge, $|a|$, and the maps $i_a,i_{\bar a}$ are both isomorphisms;
\item{$\bullet$} $Y$ consists of a single edge, $|a|$, with two distinct
vertices $\alpha(a)$, $\omega(a)$, and the images of both maps
$i_a,i_{\bar a}$ are subgroups of index 2 in the corresponding vertex groups.  

\noindent
A graph of groups over a finite graph is {\it non-elementary} if, after modifying it to a reduced graph of groups by elementary collapses, it is not
simply elementary.

\medskip
We will need the following property of non-elementary graphs of groups.

\medskip
\noindent
{\bf 4.1.6 Lemma.}
{\it Let $\cal G$ be a non-elementary graph of groups over 
a finite graph $Y$, and assume that all edge groups in $\cal G$
are finite. Let $X=X({\cal G},\Xi)$ be the Bass-Serre tree of $\cal G$.}
\item{(1)} {\it For each $v\in V_Y$ with infinite vertex group $G_v$
there is $a\in O_Y$ with $\alpha(a)=v$ such that any lift 
of $a$ to $X$ (through $\pi$) separates $X$ into two subtrees,
each of which cantains lifts of all vertices of $Y$.}
\item{(2)} {\it If all vertex groups of $\cal G$ are finite then $X$ is
an infinite locally finite tree, 
and there is $v\in V_Y$ such that any lift of $v$
to $X$ splits $X$ into at least three infinite components.}

\medskip\noindent
{\bf Proof:}
To prove part (1), fix a vertex $v\in V_Y$ for which $G_v$ is infinite.
We first claim that $v$ has more than one lift in $X$. 
If this were not the case, the unique lift $\tilde v$ of $v$ would be
fixed by all of $G=\pi_1({\cal G},\Xi)$. Hence we would have $G_v=G$,
and this could only happen if $\cal G$ was reducing to a graph of
groups over a single vertex, contradicting the assumption that $\cal G$
is non-elementary.

Now, fix two distinct lifts $v_1,v_2$ of $v$ in $X$, and let $\tilde a$
be the first oriented edge in $X$ on the unique path from $v_1$ to $v_2$.
We claim that its projection $a:=\pi(\tilde a)$ is as required.
Indeed, since any edge in $Y$ starting at $v$ lifts to infinitely many
edges in $X$ starting at $v_2$, for each vertex $u\in V_Y$ there is
its lift $\tilde u$ in $X$ such that $v_2$ lies on the path in $X$ from
$v_1$ to $\tilde u$. This shows that the subtree of $X$ obtained by
splitting at $\tilde a$ and containing $v_2$, contains also lifts
of all vertices of $Y$. The other subtree obtained by the same splitting
contains lifts of all vertices of $Y$ by a similar argument. For other
lifts of $a$ the assertion is true by transitivity of $G$ on the set of
all these lifts, and by $G$-invariance of the projection $\pi$.
This completes the proof of part (1).

To prove part (2), note that $X$ is obviously locally finite.
Moreover, $X$ is infinite since the fundamental group of any reduced
not simply elementary graph of groups is infinite.
To prove existence of a vertex $v$ as required, note that existence
of such $v$ is clearly preserved by elementary collapses. Thus,
it is sufficient to prove it in the case of reduced graphs of groups $\cal G$.
For a reduced graph of groups $\cal G$, any vertex splits the Bass-Serre
tree $X$ into as many infinite components as the valence of this vertex.
Thus, it is sufficient to show that if $\cal G$ is reduced and
non-elementary, then the Bass-Serre tree $X$ has a vertex with
valence at least 3. It is not hard to see that if $\cal G$ is reduced and
non-elementary (i.e. not simply elementary), then the underlying graph
$Y$ contains a vertex $v$ with one of the following properties:
\item{$\bullet$} there are at least two oriented edges in $Y$ starting at $v$;
\item{$\bullet$} there is an oriented edge $a$ terminating at $v$
such that the index of the subgroup $i_a(G_{|a|})<G_{v}$
is at least 3.

\noindent
In any of these two cases lifts of $v$ in $X$ have valence at least 3,
which finishes the proof.

\bigskip\noindent
{\bf 4.2. $E{\cal Z}$-structures.}

\smallskip
For completeness of the exposition, we recall the notions of 
$E{\cal Z}$-structure and $E{\cal Z}$-boundary of a group.
A slightly weaker version of this concept, callad $\cal Z$-structure,
is due to Bestvina [Be]. A generalization for groups with torsion
was introduced by Dranishnikov [Dra]. Farell and Lafont [FL]
studied an equivariant analogue, which applied only to torsion free groups.
The concept presented below generalizes all these approches,
and it has appeared in this form in Martin's paper [Ma]
(while its slightly stronger version was studied by Rosenthal [Ros]).
The concept of $E{\cal Z}$-boundary unifies and generalizes the notions
of Gromov boundary, CAT(0) boundary, and systolic boundary 
(as introduced in [OP]). Existence of an $E{\cal Z}$-structure for a group $G$
implies that $G$ satisfies the Novikov conjecture.

\medskip\noindent
{\bf 4.2.1 Definition.}
An $E{\cal Z}$-{\it structure} for a finitely generated group $G$ is a pair $(\overline{E},Z)$
of spaces (with $Z\i \overline{E}$) such that:
\item{$\bullet$} $\overline{E}$ is a Euclidean retract (i.e. a compact, contractible
and locally contractible space with finite covering dimension; such a space
is automatically metrizable);

\item{$\bullet$} $\overline{E}\setminus Z$ is a cocompact model of a 
classifying space for proper actions of $G$ (i.e. a contractible CW-complex
equipped with a properly discontinuous cocompact and cellular action of $G$,
such that for every finite subgroup $H<G$ the fixed point set 
$(\overline{E}\setminus Z)^H$ is
nonempty and contractible);

\item{$\bullet$} $Z$ is a $\cal Z$-set in $\overline{E}$ (i.e. $Z$ is a closed subspace
in $\overline{E}$ such that for any open set $U\i \overline{E}$ the inclusion $U\setminus Z\to U$
is a homotopy equivalence);

\item{$\bullet$} compact sets fade at infinity, that is, for every compact
set $K\i \overline{E}\setminus Z$, any point $z\in Z$, and any neighbourhood $U$
of $z$ in $\overline{E}$, there is a smaller neighbourhood $V\i U$ of $z$ such that
if a $G$-translate of $K$ intersects $V$ then it is contained in $U$;
\item{} this is equivalent to requiring that the set of $G$-translates of any
compact $K\i \overline{E}\setminus Z$ is a null family of subsets  in $\overline{E}$;

\item{$\bullet$} the action of $G$ on $\overline{E}\setminus Z$ extends continuously
to $\overline{E}$.

\noindent
An {\it $E{\cal Z}$-boundary} for $G$ is a space $Z$ appearing in any
$E{\cal Z}$-structure $(\overline{E},Z)$ for $G$.

\medskip
To keep track of the relationship to $G$, we will usually denote
an $E{\cal Z}$-structure for $G$ as $(\overline{EG},\partial G)$,
and the corresponding classifying space $\overline{EG}\setminus\partial G$
simply as $EG$.

\medskip
In the statement of Theorem 0.3(1) we refer also to a stronger concept
of the boundary, as defined below.

\medskip\noindent
{\bf 4.2.2 Definition}.
An $E{\cal Z}$-structure $(\overline{E},Z)$ {}
is {\it strong in the sense of Carl\-sson-Pedersen} if  for each finite subgroup $H<G$
the fixed point set $Z^H$ is either empty or a $\cal Z$-set in $\overline{E}^H$.
An {\it $E{\cal Z}$-boundary strong in the sense of Carlsson-Pedersen} for $G$ is a space $Z$ appearing in any
$E{\cal Z}$-structure $(\overline{E},Z)$ for $G$ strong in the sense of Carlsson-Pedersen.

\medskip
The above concept strengthens slightly, in a natural way, the concepts
appearing in the works of Carlsson and Pedersen [CP], as well as
Rosenthal [Ros]. It has appeared in Martin's paper [Ma],
where it turned out to be natural from the point of view of the
combination theorem being the main result of that paper.


\bigskip\noindent
{\bf 4.3. An $E{\cal Z}$-structure for a graph of groups with finite
edge groups.}

\smallskip
Let $\cal G$ be a graph of groups as in Theorem 0.3(1),
over a finite graph $Y$.
It means that all edge groups $G_e$ in $\cal G$ are finite, and each
vertex group $G_v$ is equipped with an $E{\cal Z}$-structure
$(\overline{EG_v},\partial G_v)$. In this subsection we briefly recall
the construction
of an $E{\cal Z}$-structure
$(\overline{E_MG},\partial_M G)$ for the fundamental group 
$G=\pi_1({\cal G},\Xi)$. This is a rather special case of a much more
general construction presented by Alexandre Martin in 
[Ma]. Our description
is adapted to the case under our interest.

Apart from the $E{\cal Z}$-structures $(\overline{EG_v},\partial G_v)$,
as initial data for the construction we need the following:
for each oriented edge $a\in O_Y$ we choose a point $p_a\in EG_{\omega(a)}$
which is fixed by the subgroup $i_a(G_{|a|})<G_{\omega(a)}$. Note that
the subgroup $i_a(G_{|a|})$ is finite, and hence its fixed point set
in $EG_{\omega(a)}$ is not empty, which justifies existence of $p_a$.
The tuple of data $( \{ (\overline{EG_v},\partial G_v) :v\in V_Y \}, 
\{ p_a:a\in O_Y \} )$
as above is an example of {\it an $E{\cal Z}$-complex of classifying spaces compatible with
$\cal G$}, see Definitions 2.2 and 2.6 in [Ma].

\smallskip
We first describe a cocompact model $E_MG$ of a classifying space for
proper actions of $G$. As in Remark 4.1.4, 
for each $a\in O_Y$, let $\tau_{|a|}$ be
a nonoriented edge with its two associated oriented edges $\tau_a$
and $\overline{\tau_a}=\tau_{\bar a}$, and suppose that its endpoints
$\alpha(\tau_a)$ and $\omega(\tau_a)$ are distinct. View $\tau_{|a|}$
as a topological space homeomorphic to a segment. Put
$$
E_MG:=\big[ G\times\big(  (\bigsqcup_{v\in V_Y}EG_v)\sqcup
(\bigsqcup_{e\in|O|_Y}\tau_e) \big) \big] /\sim,
$$
where the equivalence relation $\sim$ is induced by the following equivalences:
\item{$\bullet$} $(gh,x)\sim(g,hx)$ for all $g\in G$, $v\in V_Y$, $x\in EG_v$
and $h\in G_v$;
\item{$\bullet$} $(gh,y)\sim(g,y)$ for all $g\in G$, $e\in|O|_Y$, $y\in\tau_e$
and $h\in G_e$;
\item{$\bullet$}  $(g,p_a)\sim(g,\omega(a))\in G\times\tau_{|a|}$
for all $g\in G$ and all $a\in O_Y:a^+\subset \Xi$;
\item{$\bullet$}  $(gs_a^{-1},p_a)\sim(g,\omega(a))\in G\times\tau_{|a|}$
for all $g\in G$ and all $a\in O_Y:a^+\not\subset \Xi$.

\noindent
The action of $G$ on $E_MG$ is induced by $h\cdot(g,x)=(hg,x)$
for any 
$x\in (\bigsqcup_{v\in V_Y}EG_v)\sqcup
(\bigsqcup_{e\in|O|_Y}\tau_e)$ and any $g,h\in G$.
This is a specification of the construction from Section II.2 in [Ma].
Theorem II.2.3 in the same paper asserts that $E_MG$ is indeed
a cocompact model of a classifying space for proper actions of $G$
(which also can be easily seen directly in this rather special case).

In addition to the above, we have a 
continuous 
$G$-equivariant map
$p:E_MG\to X({\cal G},\Xi)$ to the Bass-Serre tree, induced by
$p((g,x))=(gG_v,v)$ for $x\in EG_v$, and by $p((g,y))=(gG_{|a|},y)
\in\{ (gG_{|a|} \}\times\tau_{|a|}$
for $y\in\tau_{|e|}$ (where in the last expression we refer
to the description of $X=X({\cal G},\Xi)$ as in Remark 4.1.4).

Note that for each vertex $t\in V_X$ the preimage $p^{-1}(t)$
is a subspace of $E_MG$ which is an embedded copy of $EG_{\pi(t)}$.
This subspace will be denoted $EG_t$, which nicely interplays with
the following. If we denote by $G_t$ the subgroup of $G$ stabilizing
the vertex $t$, then $EG_t$ is invariant under $G_t$, and it is
a classifying space for proper actions for $G_t$.

Note also that for an edge $\varepsilon=[gG_{|a|},\tau_{|a|}]$ of $X$,
denoting by $\varepsilon^\circ$ its geometric interior,
the closure in $E_MG$ of the preimage $p^{-1}(\varepsilon^\circ)$,
denoted $\overline{p^{-1}(\varepsilon^\circ)}$, 
is an embedded copy of $\tau_{|a|}$.
We call each set of this form {\it a segment} in $E_MG$.
The endpoint of this segment, which projects through $p$ to $t=[gG_{|a|},\omega(\tau_a)]$,
belongs to the subspace $EG_t$, and we call it {\it the attaching point
of the segment $\overline{p^{-1}(\varepsilon^\circ)}$ 
in the subspace $EG_t$}. 
Observe also that $p$ establishes a bijective correspondence between
the nonoriented edges of $X$ and the segments of form 
$\overline{p^{-1}(\varepsilon^\circ)}$ 
as above. We will call the segment 
$\overline{p^{-1}(\varepsilon^\circ)}$ in $E_MG$ {\it the lift of the edge $\varepsilon$}
of $X$.

\smallskip
We now pass to the description of a set $\partial_{{Stab}}G$,
which is a part of $\partial_MG$. This is the specialization to our case
of the construction given at the end of Section 2.1 in [Ma].
Put
$$
\partial_{Stab}G:=\big(  G\times (\bigsqcup_{v\in V_Y}\partial G_v) \big) / \sim,
$$
where $\sim$ is induced by the equivalences
$(gh,x)\sim(g,hx)$ for all $g\in G$, $v\in V_Y$, 
$x\in\partial G_v$ and $h\in G_v$. 
The action of $G$ on $\partial_{Stab}G$ is given by acting from the left on
the first coordinate. 
We also have the $G$-equivariant projection
$p_{Stab}:\partial_{Stab}G\to V_X$ induced by $p_{Stab}((g,x))=(g,v)$
for all $v\in V_Y$ and all $x\in\partial G_v$.
For any vertex $t\in V_X$, the preimage $p_{Stab}^{-1}(t)$ is $G_t$-invariant
and has a (unique up to $G_t$-action) identification with the boundary
$\partial G_{\pi(t)}$. We denote this preimage by $\partial G_t$.
The union $\overline{EG_t}:=EG_t\sqcup\partial G_t$ 
has a (unique up to $G_t$-action)
identification with $\overline{EG_{\pi(t)}}$. Under the topology induced
from this identification, the pair
$(\overline{EG_{t}},\partial G_t)$ is an $E{\cal Z}$-structure for $G_t$.

\smallskip
A third ingredient in the description of the $E{\cal Z}$-structure
$(\overline{E_MG},\partial_MG)$ is the set $\partial X$
of ends of the Bass-Serre tree $X=X({\cal G},\Xi)$. More precisely,
this is the set of equivalence classes of infinite combinatorial rays in $X$
for the relation of coincidence except at possibly some finite initial parts.
The action of $G$ on $X$ induces the action on $\partial X$.
We then put
$$
\partial_MG:=\partial_{Stab}G\sqcup\partial X \hbox{\quad and \quad}
\overline{E_MG}:=E_MG\sqcup\partial_MG.
$$
The union of the maps $p$, $p_{Stab}$ and the identity map on $\partial X$
gives the map $\bar p:\overline{E_MG}\to X\sqcup\partial X$ which is $G$-equivariant.
Moreover, for each vertex $t\in V_X$, the preimage $(\bar p)^{-1}(t)$
coincides with $\overline{EG}_t$.

\medskip
We now recall the topology in $\overline{E_MG}$, as described 
in Section IV.5 of
[Ma]. In fact, we are interested only in the restricted topology
in the boundary $\partial_MG$, so we recall only this part of the information.
We do this by describing, for any point $z\in\partial_M G$, a basis of open
neighbourhoods of $z$ in $\partial_M G$.

Fix a vertex $t_0$ in the Bass-Serre tree $X$.
If $z\in\partial X$, for any integer $n\ge1$ let $X_n(z)$ be the subtree
of $X$ spanned on all vertices $t\in V_X$ for which the path in $X$
connecting $t_0$ to $t$ has the same first $n$ edges as 
the infinite path in $X$ from $t_0$ to $z$. Denote by $\partial X_n(z)$
the set of ends in this subtree, viewing it canonically as
a subset of $\partial X$. Put 
$V_n(z):=p_{Stab}^{-1}(V_{X_n(z)})\cup\partial X_n(z)$.
As a basis of open neighbourhoods of $z$ in $\partial G$ take the
family of sets $V_n(z)$ for all integer $n\ge1$.

If $z\in\partial_{Stab}G$, let $t$ be the vertex of $X$ such that 
$z\in\partial G_t$. Let $U$ be an open neighbourhood of $z$ in 
$\overline{EG_t}$ (for the topology induced from the identification
with $\overline{EG_{\pi(t)}}$). Put $\widetilde V_U$ to be the set
of all elements $u\in\partial G$ with $p(u)\ne t$ and such that
the geodesic in $X\cup\partial X$ from $t$ to $p(u)$ starts with an edge $\varepsilon$
which lifts through $p$ to a segment in $E_MG$ whose attaching point
in $EG_t$ belongs to $U$. Put then $V_U(z):=U\cup\widetilde V_U$.
As a basis of open neighbourhoods of $z$ in $\partial G$ take the
family of sets $V_U(z)$, where $U$ runs through some basis of open
neighbourhoods of $z$ in $\overline{EG_t}$.

\bigskip\noindent
{\bf 4.4. Proof of Theorem 0.3(1).}

\medskip
Part (1) of Theorem 0.3 is a direct consequence of the following
property of $E{\cal Z}$-boundaries $\partial_MG$ described in
the previous subsection.

\medskip\noindent
{\bf 4.4.1 Lemma.}
{\it Under assumptions of Theorem 0.3(1), we have 
$$\partial_MG\cong\widetilde\sqcup(\partial G_{v_1},\dots,\partial G_{v_k}).$$}

\noindent
{\bf Proof:}
Let $\cal G$ be a graph of groups as in Theorem 0.3(1),
and let $\partial_MG$ be the $E{\cal Z}$-boundary
of the fundamental group $G=\pi_1({\cal G},\Xi)$,
as described in the previous subsection. 
Consider first the special case when all vertex groups $G_{v_i}$
are finite. It follows from the definition of an $E{\cal Z}$-structure
that the boundaries $\partial G_{v_i}$ are then all empty.
Thus, by our convention, we have that the amalgam
$\widetilde\sqcup(\partial G_{v_1},\dots,\partial G_{v_k})$
is then the Cantor space $C$. On the other hand, the boundary
$\partial_MG$ reduces in
this case to the part $\partial X$. By Lemma 4.1.6(2), $X$ is then
an infinite uniformly locally finite tree 
such that the set of
vertices splitting it into at least three infinite components
is a net in $X$ (i.e. there is $D>0$ such that every vertex of $X$
remains at combinatorial distance at most $D$ from a vertex
in this set).
A straightforward argument shows that $\partial X$,
with the topology described in the previous subsection,
is then homeomorphic to the Cantor space $C$.
Thus the theorem follows in the considered case. 

We now pass to the case when at least one vertex group is infinite.
Recall that, by definition, an $E{\cal Z}$-boundary of a group
(if exists) is nonempty iff the group is infinite.
Without loss of generality, suppose that for some $m\in\{ 1,\dots,k\}$
the vertex groups $G_{v_1},\dots,G_{v_m}$ are infinite,
while the remaining ones are finite.
Since the boundaries $\partial G_{v_j}$ for $j>m$ are empty,
by our convention we have
$\widetilde\sqcup(\partial G_{v_1},\dots,\partial G_{v_m},\dots,
\partial G_{v_k})=
\widetilde\sqcup(\partial G_{v_1},\dots,\partial G_{v_m})$.
Thus we need to show that 
$\partial_MG$ is homeomorphic to 
$\widetilde\sqcup(\partial G_{v_1},\dots,\partial G_{v_m})$.

By definition of $E{\cal Z}$-boundary (Definition 4.2.1),
$\partial_MG$ is compact and metrizable.
Using the notation
introduced in the previous subsection,
define a family ${\cal Y}={\cal Y}_1\sqcup\dots\sqcup{\cal Y}_m$
of subsets in $\partial_M G$ as follows.
For each $i\in \{ 1,2,\dots,m \}$ put
${\cal Y}_i:=\{ \partial G_t:t\in V_X, \pi(t)=v_i \}$.
In view of Theorem 0.2, it is sufficient to show that $\cal Y$
is a $(\partial G_{v_1},\dots,\partial G_{v_m})$-regularizing
family for $\partial_MG$. Thus we need to check conditions
(a1)-(a5) of the introduction.

\medskip
Recall that the topology in the subspace $\partial G_t$ 
induced from that in $\partial_MG$ coincides with the
topology provided by the identification of $\partial G_t$
with $\partial G_{\pi(t)}$
(compare Proposition IV.5.19 in [Ma]). 
Thus, each subset $\partial G_t$
is an embedded copy of $\partial G_{\pi(t)}$, which
verifies condition (a1).

To check condition (a2), i.e. nullnes of the family $\cal Y$,
we need to show that for each open covering $\cal U$ of 
$\partial_MG$ there is a finite subfamily ${\cal A}\subset{\cal Y}$
such that for every $Z\in{\cal Y}\setminus{\cal A}$ there is $U\in{\cal U}$
that contains $Z$. Obviously, without loss of generality we may
assume that $\cal U$ is finite and consists of sets from bases of open 
neighbourhoods of points. Suppose that
$$
{\cal U}=\{ V_{U_1}(z_1),\dots,V_{U_p}(z_p),V_{n_1}(z_1'),\dots,
V_{n_q}(z_q')\}.
$$
It is not hard to see that the above property holds for
$\cal U$ with ${\cal A}=\{ \partial G_{p(z_1)},\dots,\partial G_{p(z_p)} \}$.
We omit further details.

To check (a3), choose any $Z\in{\cal Y}$, i.e. a subset $\partial G_t$
for some vertex $t$ of $X$ such that $\pi(t)=v_i$ and $i\le m$.
Choose also any point $z\in\partial G_t$, any open neighbourhood $U$
of $z$ in $\overline{EG_t}$, and consider the associated open
neighbourhood $V_U(z)$ from the local basis at $z$ in $\partial_MG$,
as described at the end of Subsection 4.3. We need to show that
$V_U(z)$ contains a point of $\partial_MG\setminus\partial G_t$.

Recall that we denote by $G_t$ the subgroup of $G$ stabilizing $t$, 
and that this sungroup is isomorphic to $G_{v_i}$, and hence it is
infinite. Moreover, the pair $(\overline{EG_t},\partial G_t)$ is
an $E{\cal Z}$-structure for $G_t$.
Since $\cal G$ is non-elementary, it follows from Lemma 4.1.6(2)
that some edge $\varepsilon$ of $X$ adjacent to $t$ splits $X$ into subtrees
containing lifts of all vertices of $Y$.
Let $\tau=\overline{p^{-1}(\varepsilon^\circ)}$ be the segment
of $E_MG$ which is the lift of $\varepsilon$, and let $x$ be the
attaching point of $\tau$ in $EG_t$.
Since, by definition of an $E{\cal Z}$-structure, compact subsets
of $EG_t$ fade at infinity, we get that there is $x'$ in $G_t$-orbit of $x$
such that $x'\in U$. This $x'$ is the attaching point in $EG_t$
of another segment $\tau'$ of $E_MG$. By $G$-equivariance,
the image $\varepsilon':=p(\tau')$ is a different from $\varepsilon$
edge of $X$ adjacent to $t$ that splits $X$ into subtrees containing
lifts of all vertices of $Y$.
Let $s$ be a vertex of $X$ which is a lift of $v_1$, and which after splitting
$X$ at $\varepsilon'$ belongs to the other component than $t$.
By definition of $V_U(z)$, we see that $\partial G_s\i V_U(z)$.
Since $\partial G_s\ne\emptyset$ (because $\partial G_s\cong
\partial G_{v_1}$, and $G_{v_1}$ is infinite), this completes
the verification of (a3).

The argument in the previous paragraph shows in fact that for
each $i\in\{ 1,\dots,m \}$ any point of $\partial_{Stab}G$ belongs
to the closure in $\partial_MG$ of the subset
$\cup{\cal Y}_i=p_{Stab}^{-1}(\pi^{-1}(v_i))$.
To check condition (a4), i.e. that $\cup{\cal Y}_i$ is dense in 
$\partial_MG$, it remains to show that any point of $\partial X$
also belongs to the closure of $\cup{\cal Y}_i$.
Let $z\in\partial X$, and let $V_n(z)$ be a neigbourhood of $z$ in
$\partial_MG$ which belongs to a local basis at $z$, as described at
the end of Subsection 4.3. Let $D$ be the combinatorial diameter of
the graph $Y$, and let $u$ be the vertex on the infinite path in $X$
from $t_0$ to $z$, at distance $n+D$ from $t_0$. Let $s$ be a vertex
of $X$ which is a lift of $v_i$ lying at combinatorial dostance $\le D$
from u (by definition of $D$, such $s$ always exists).
Observe that $s\in V_{X_n(z)}$, and hence 
$\partial G_s=p_{Stab}^{-1}(s)\i V_n(z)$.
Since $\partial G_s\ne\emptyset$ (because $\partial G_s\cong
\partial G_{v_i}$), this completes the verification
of condition (a4).

To check condition (a5), we make the following two observations,
the direct proofs of which we omit. First, note that for any 
$z\in\partial X$, any set $V_n(z)$ from the local basis at $z$ is both open 
and closed in $\partial_MG$. Second, observe that any two points
of $\partial_MG$ not contained in the same set $Z\in{\cal Y}$
(i.e. in the same set $\partial G_t$ for any 
$t\in p^{-1}(\{ v_1,\dots,v_m \})$) can be separated from each other
by some set $V_n(z)$, for appropriately chosen $z$ and $n$. This completes the verification
of condition (a5), and thus completes the proof.

\bigskip
\noindent
{\bf 5. Gromov boundaries and $\hbox{CAT}(0)$ boundaries.} 

\medskip
In this section we prove parts (2) and (3) of Theorem 0.3. 
It is not hard to give direct proofs of these results,
by referring to the characterization of dense amalgams provided
in Theorem 0.2. However, we present shorter arguments, based
on properties of $E{\cal Z}$-boundaries $\partial_MG$
constructed in Subsection 4.3.

\medskip\noindent
{\it 
Gromov boundary and the proof of Theorem 0.3(2).}

\smallskip

We use the following result of A. Martin (see Corollary 9.19 in [Ma]).

\medskip\noindent
{\bf 5.1 Lemma.}
{\it Let $\cal G$ be a graph of groups satisfying the assumptions
of part (2) of Theorem 0.3.
Let $(\overline{PG_{v_i}},\partial G_{v_i})$
be the $E{\cal Z}$-structures for the vertex groups $G_{v_i}$
provided by the compactifications of appropriate Rips complexes
$PG_{v_i}$
by means of Gromov boundaries $\partial G_{v_i}$ of these groups.
Then the $E{\cal Z}$-boundary $\partial_MG$ for $G=\pi_1({\cal G})$
obtained from the above $E{\cal Z}$-structures as in Subsection 4.3
is $G$-equivariantly homeomorphic to the Gromov boundary of $G$.}

\medskip
Note that, under assumptions of the
above lemma, it follows from Lemma 4.4.1 that
$\partial_MG\cong
\widetilde\sqcup(\partial G_{v_1},\dots,\partial G_{v_k})$.
Consequently,
Theorem 0.3(2) follows from Lemma 5.1.

\bigskip\noindent
{\it
$\hbox{CAT}(0)$ boundary and the proof of Theorem 0.3(3).}

\smallskip
Recall that if a group $\Gamma$ acts geometrically (i.e. by isometries, properly discontinuously and cocompactly) on a $\hbox{CAT}(0)$
space $W$, and if $\overline W$ denotes the compactification of $W$
by means of its $\hbox{CAT}(0)$ boundary $\partial W$, then the pair
$(\overline{W},\partial W)$ is an $E{\cal Z}$-structure for $\Gamma$.

We work under assumptions and notation of Theorem 0.3(3). 
Let $(\overline{E_MG},\partial_MG)$ be the $E{\cal Z}$-structure
for $G=\pi_1(\cal G)$ constructed as in Subsection 4.3 out of
$\hbox{CAT}(0)$ $E{\cal Z}$-structures $(\overline{\Delta_i},\partial\Delta_i)$.
We make the following observations concerning this $E{\cal Z}$-structure.

\medskip\noindent
{\bf 5.2 Lemma.}
\item{(1)}
{\it The space $E_MG$ carries a natural geodesic metric for which
it is $\hbox{CAT}(0)$, and for which $G$ acts on $E_MG$ geometrically.}
\item{(2)} 
{\it The boundary $\partial_MG$ naturally coincides (as a topological space)
with the $\hbox{CAT}(0)$ boundary $\partial E_MG$ (for the $\hbox{CAT}(0)$
geodesic metric 
in $E_MG$ as in part (1)).}

\medskip\noindent
{\bf Proof:}
To prove (1), note that $E_MG$ is obtained from copies of the 
$\hbox{CAT}(0)$
spaces $\Delta_i$, and from copies of the segment, by gluing the endpoints 
of the segments to the appropriate attaching points in copies of $\Delta_i$.
By putting at each segment the standard euclidean metric of length 1,
we get on $E_MG$ the induced length metric which is geodesic
(see I.5.26 in [BH]). Since we perform the gluings along singletons,
which are obviously convex as subspaces, the successive application
of Basic Gluing Theorem 11.1 of [BH] shows that $E_MG$ with
the above metric is $\hbox{CAT}(0)$. Obviously, with this metric $G$ acts on
$E_MG$ by isometries. The action is proper and cocompact by
definition of $E{\cal Z}$-structure.

To prove part (2), choose a base point $x_0\in E_MG$ as a point
in some copy of some $\Delta_i$. There are two kinds of geodesic rays
in $E_MG$ starting at $x_0$:

\item{(a)} those which pass through infinitely many segments;
\item{(b)} those which, after passing through finitely many segments,
eventually coincide with a geodesic ray in some copy of some $\Delta_i$.

\noindent
We define a map $h:\partial E_MG\to\partial_MG$ as follows.
If $\xi\in\partial E_MG$ is represented by a geodesic ray of kind (a)
above, note that the sequence of segments through which this ray
successively passes lifts to a sequence of edges in the Bass-Serre
tree $X$ which forms an infinite combinatorial ray $\rho$;
denoting by $[\rho]\in\partial X$ the end of $X$ represented by $\rho$,
we put $h(\xi):=[\rho]$. If $\xi$ is represented by a geodesic ray
of kind (b), its final part (which is a geodesic ray in some copy $EG_t$
of some $\Delta_i$) induces a point $z$ in the $\hbox{CAT}(0)$ boundary of this copy
(i.e. a point in $\partial G_t\subset\partial_{Stab}G$); we then put
$h(\xi):=z$. The so described map 
$h:\partial E_MG\to \partial_MG=\partial_{Stab}G\cup\partial X$ 
is easily seen to be a bijection.
As both spaces $\partial E_MG$ and $\partial_MG$ are compact,
to finish the proof of (2) we need to show that $h$ is continuous.

Recall (e.g. from II.8.6 in [BH]) that a point $\xi$ of the
boundary of a $\hbox{CAT}(0)$ space $W$, represented by a geodesic ray
$\gamma_\xi$ started at a point $x_0\in W$, has a basis of open
neighbourhoods of form 
$$
{\cal U}(\gamma_\xi,r,\varepsilon)=\{ \eta\in\partial W : 
d_W(\gamma_\xi(r),\gamma_{\eta}(r))<\varepsilon \},
$$
where $r$ and $\varepsilon$ run through arbitrary  positive real numbers,
$\gamma_\eta$ is the geodesic ray in $W$ started at $x_0$ and
representing $\eta$, and $\gamma_\eta(r)$ is the point on $\gamma_\eta$ 
at distance
$r$ from $x_0$.
Below we will make use of the sets of the above form 
${\cal U}(\gamma_\xi,r,\varepsilon)$ for the space $W=E_MG$.
Without loss of generality, we assume that the base point 
$x_0\in E_MG$
is chosen in the subspace $EG_{t_0}$ (which is a copy of some $\Delta_i$), where $t_0$ is the base vertex
in the Bass-Serre tree $X$, as
fixed at the end of Subsection 4.3 (in the description of local bases of neighbourhoods for the topology in $\partial_MG$).
Let $p=h(\xi)$ be any point of $\partial_MG$, and $V$ its any open neighbourhood. We need to indicate an open neighbourhood 
$\cal U$ of $\xi$
in $\partial E_MG$ such that $h({\cal U})\i V$. Clearly, we may restrict ourselves
to the case when $V$ belongs to the basis of local neighbourhoods at $p$,
as described at the end of Subsection 4.3. 

We consider two cases.
First, suppose that $p=[\rho]\in\partial X\i\partial_MG$, where $\rho$
is the combinatorial ray in $X$ induced by a geodesic ray $\gamma$
in $E_MG$ of kind (a) above; then $\gamma$ starts at $x_0$ and represents $\xi$. Let $V=V_n(\xi)$ for some $n\ge1$. 
Let $x_n$ be the most distant from $x_0$ point on the $n$-th
segment in $E_MG$ traversed by the ray $\gamma$, and let
$r=d_{E_MG}(x_0,x_n)$. It is then easy to see that
${\cal U}={\cal U}(\gamma, r,1)$ is as required. We omit further details.
In the second case, suppose that $p=z\in\partial G_t$ for some vertex
$t$ of the Bass-Serre tree $X$. Then $\xi$ is represented by the
geodesic ray $\gamma$ in $E_MG$ started at $x_0$, which eventually
coincides with the geodesic ray $\gamma_t$ in $EG_t$ representing $z$ and started at the attaching point
$x_t$ of the segment through which any geodesic ray started at $x_0$ 
enters $EG_t$. Let $r_t=d_{E_MG}(x_0,x_t)$. 
Let also $V=V_U(z)$ for some open neinghbourhood $U$ of $z$
in $\overline{EG_t}$. By the description of the topology in the $\hbox{CAT}(0)$ compactification
$\overline{EG_t}$ (see again II.8.6 in [BH]), there are positive reals
$r$ and $\varepsilon$ with the following property: for any geodesic ray
$\beta$ in $EG_t$ started at $x_t$, if 
$d_{EG_t}(\beta(r),\gamma_t(r))<\varepsilon$ then for any
$r'\in(r,\infty]$ the point $\beta(r')$ belongs to $U$ (here,
by $\beta(\infty)$ we mean the point in the boundary represented by
$\beta$). It is not hard to see that then the ball of radius $\varepsilon$
in $EG_t$ centered at 
$\gamma_t(r+\varepsilon)=\gamma(r_t+r+\varepsilon)$ is also contained
in $U$. From this, it follows fairly directly that the set
${\cal U}={\cal U}(\gamma,r_t+r+\varepsilon,\varepsilon)$
is as required. This completes the proof of the lemma.

\medskip
Now, since by Lemma 4.4.1 under our assumptions
we have 
$\partial_MG\cong\widetilde\sqcup(\partial\Delta_1,\dots,\partial\Delta_k)$, Theorem 0.3(3) follows from Lemma 5.2
by putting $\Delta=E_MG$.

\bigskip
\noindent
{\bf 6. Systolic boundaries.} 

\medskip
In this section we prove part (4) of Theorem 0.3.
In Subsection 6.1 we briefly recall the definition and basic properties
of systolic complexes and groups. In Subsection 6.2 we construct
the systolic complex $\Sigma$ appearing in the assertion
of Theorem 0.3(3), as appropriate tree of systolic complexes.
In Subsection 6.3 we recall the concept of systolic boundary.
Finally, in Subsection 6.4 we prove Theorem 0.3(4) by studying
the systolic boundary of the earlier described complex $\Sigma$
in the light of the characterization of dense amalgams provided in Theorem 0.2. 

\bigskip\noindent
{\bf 6.1 Systolic complexes and groups.}

\smallskip
Systolic complexes have been introduced in the
paper by T.~Ja\-nusz\-kiewicz and the author [JS]. These are
the simply connected simplicial complexes of arbitrary 
dimension that satisfy some local (combinatorial) condition
that resembles nonpositive curvature. 
A group is called {\it systolic} if it acts geometrically
(i.e. by simplicial automorphisms, properly discontinuously
and cocompactly) on a systolic complex.
It is shown in [JS] that systolic groups are
biautomatic, and hence also semihyperbolic, and that they appear
in abundance in arbitrary (virtual) cohomological dimension. 

We recall briefly the definition of a systolic complex.
A simplicial complex is {\it flag} if its any set of vertices pairwise
connected with edges spans a simplex.
A {\it full cycle} in a simplicial complex is a full subcomplex
isomorphic to a triangulation of the circle $S^1$.
A simplicial complex is {\it $6$-large} if it is flag
and contains no full cycle with less than $6$ edges.
A simplicial complex is {\it systolic}
if it is simply connected and its link at any vertex 
is $6$-large. This simple definition describes spaces with
surprisingly rich geometric (but expressed in purely
combinatorial terms) structure. One of the basic observations is that
systolic complexes are contractible, which is an analogue of
Cartan-Hadamard theorem. We recall few further facts that we 
need in the present paper.
The first result below is due to Victor Chepoi and Damian Osajda.

\medskip\noindent
{\bf 6.1.1 Theorem} (Theorem C in [ChO]).
{\it Let $H$ be a finite group acting by automorphisms on a
locally finite systolic complex $\Upsilon$. Then $\Upsilon$
contains a simplex which is $H$-invariant.}

\medskip
The next result concerns existence of natural $E{\cal Z}$-structures
for systolic groups.

\medskip\noindent
{\bf 6.1.2 Theorem} (Theorem A in [OP] and Theorem E in [ChO]).
{\it Let $\Upsilon$ be a systolic complex acted upon geometrically
by a systolic group $G$. Then there is a compactification
$\overline\Upsilon=\Upsilon\cup\partial\Upsilon$ such that the pair
$(\Upsilon,\partial\Upsilon)$ is an $E{\cal Z}$-structure for $G$.}

\medskip
The paper [OP] by Damian Osajda and Piotr Przytycki contains
construction of a compactification $\overline\Upsilon$ as above,
and the corresponding space 
$\partial\Upsilon=\overline\Upsilon\setminus\Upsilon$ resulting
from this construction is {\it the systolic boundary} of $\Upsilon$,
as appearing in the statement of Theorem 0.3(4).
It is shown in [OP] that if a systolic group $G$ is word hyperbolic
then its any systolic boundary (depending on the choice of
a systolic complex on which $G$ acts geometrically) coincides
with the Gromov boundary of $G$.
Thus, this notion naturally extends the concept of ideal boundary
to the class of systolic groups which are not word hyperbolic.
In Subsection 6.3 we indicate the features of systolic boundaries
necessary for the proof of Theorem 0.3(4), and the proof itself is provided
in Subsection 6.4.

\medskip
Recall from [P] the following natural notion and a related observation.

\medskip\noindent
{\bf 6.1.3 Definition.}
A {\it tree of systolic complexes} is a simplicial complex $\Upsilon$
equipped with a simplicial map $p:\Upsilon\to X$ onto a simplicial
tree $X$ satisfying the following. For every vertex $t$ of $X$
the preimege $p^{-1}(t)$ is a systolic complex, and for every
open edge $e^\circ$ of $X$  the closure in $\Upsilon$ of the preimage
$p^{-1}(e^\circ)$ is a simplex.

\medskip\noindent
{\bf 6.1.4 Lemma} ([P], Section 7).
{\it If $p:\Upsilon\to X$ is a tree of systolic complexes then
$\Upsilon$ is itself a systolic complex.}

\bigskip\noindent
{\bf 6.2 Graphs of systolic groups and the construction
of $\Sigma=\Sigma({\cal G},\Xi)$.}

\medskip
In this subsection, under assumptions of Theorem 0.3(4),
we construct a systolic complex $\Sigma$ as asserted in the
theorem. Verification that the systolic boundary $\partial\Sigma$
is homeomorphic to the appropriate dense amalgam
will be provided in Subsection 6.4.

Let $\cal G$ be as in the assumptions of Theoerem 0.3(4),
and let $\Xi$ be a maximal tree in the first barycentric subdivision $Y'$
of the underlying graph $Y$ of $\cal G$. We use the notation as in Subsection 4.1
concerning $\cal G$ and the associated objects.  

For each $a\in O_Y$, fix an embedding $j_a:\sigma_a\to\Sigma_{\omega(a)}$ of an abstract simplex $\sigma_a$
onto some simplex of $\Sigma_{\omega(a)}$ preserved by the restricted
action on $\Sigma_{\omega(a)}$ of the subgroup 
$i_a(G_{|a|})<G_{\omega(a)}$. Since the subgroup $i_a(G_{|a|})$
is finite, existence of such an embedding is ensured by
Theorem 6.1.1. We denote by $j_a^{-1}$ the inverse isomorphism
from the simplex $j_a(\sigma_a)$ to $\sigma_a$.
For each $a\in O_Y$, put $\kappa_{|a|}:=\sigma_a*\sigma_{\bar a}$
(i.e. the simplicial join of the simplices $\sigma_a,\sigma_{\bar a}$).
Consider the action of the group $G_{|a|}$ on the simplex $\kappa_{|a|}$
by simplicial automorphisms defined on vertices $u$ by
$g\cdot u:=j_a^{-1}(i_a(g)\cdot j_a(u))$
for $u\in\sigma_a$, and
$g\cdot u:=j_{\bar a}^{-1}(i_{\bar a}(g)\cdot j_{\bar a}(u))$
for $u\in\sigma_{\bar a}$. Put
$$
\Sigma=\Sigma({\cal G},\Xi):=
[G\times\big( (\bigsqcup_{v\in V_Y}\Sigma_v)\sqcup 
(\bigsqcup_{e\in|O|_Y}\kappa_e) \big)]/\sim
$$
where the equivalence relation $\sim$ is induced by the following
equivalences:
\item{$\bullet$} $(gh,x)\sim(g,hx)$ for all $g\in G$, $v\in V_Y$, 
$x\in \Sigma_v$ and $h\in G_v$;
\item{$\bullet$} $(gh,y)\sim(g,hy)$ for all $g\in G$, $e\in|O|_Y$,
$y\in\kappa_e$ and $h\in G_e$;
\item{$\bullet$} $(g,j_a(y))\sim(g,y)$
for all $g\in G$, all $a\in O_Y \hbox{ such that }a^+\subset \Xi$, and all $y\in\sigma_a\subset\kappa_{|a|}$;
\item{$\bullet$} $(gs_a^{-1},j_a(y))\sim(g,y)$
for all $g\in G$, all $a\in O_Y \hbox{ such that }a^+\not\subset \Xi$, and all $y\in\sigma_a\subset\kappa_{|a|}$.

\noindent
We denote by $[g,x]$ the equivalence class under the relation $\sim$
of an element $(g,x)$.

Observe that the above described space $\Sigma$ carries a natural
induced structure of a simplicial complex. More precisely, the injective
images in $\Sigma$ (through the quotient map provided by $\sim$)
of the simplices $\{ g \}\times\sigma$ in the copies $\{ g \}\times\Sigma_v$
or $\{ g \}\times\kappa_e$ yield the structure of a simplicial complex
for $\Sigma$. We denote the image simplices as above by $[g,\sigma]$.

$\Sigma$ comes equipped with a simplicial projection map
$p:\Sigma\to X=X({\cal G},\Xi)$ onto the Bass-Serre tree of $\cal G$.
This map is determined by its restriction to vertices, which is described
as follows: $p([g,w])=(gG_v,v)$ for any $v\in V_Y$, any 
vertex $w\in\Sigma_v$, and any $g\in G$.
For any vertex $t=(gG_v,v)$ of $X$, the preimage $p^{-1}(t)$
is a subcomplex in $\Sigma$ isomorphic to $\Sigma_v$,
and we denote it $\Sigma_t$. Similarly, for any geometric edge 
$\varepsilon=[gG_{|a|},\tau_{|a|}]$ of $X$, closure in $\Sigma$ of 
the preimage $p^{-1}(\varepsilon^\circ)$ of its interior $\varepsilon^\circ$
is a simplex of $\Sigma$, naturally isomorphic with the simplex 
$\kappa_{|a|}$. We denote this simplex by $\kappa_\varepsilon$.
As a consequence, $p:\Sigma\to X$ is a tree of systolic complexes,
as in Definition 6.1.3, and hence, by Lemma 6.1.4, $\Sigma$ is
a systolic complex.

Consider the simplicial action of the fundamental group 
$G=\pi_1({\cal G},\Xi)$ on $\Sigma$ which is described on vertices by
$h\cdot[g,w]=[hg,w]$ for any $v\in V_Y$, any vertex $w\in\Sigma_v$,
and any $g,h\in G$. This action is easily seen to be cocompact,
as it is  not hard to indicate a finite set of representatives of orbits
for the induced action of $G$ on the set of all simplices of $\Sigma$.
Moreover, the stabilizer of a vertex $[g,w]$ of $\Sigma$,
where $w$ is a vertex of $\Sigma_v$ for some $v\in V_Y$,
coincides with the subgroup 
$g\cdot\hbox{Stab}_{G_v}(w)\cdot g^{-1}<G$
(under the natural interpretation of $G_v$ as a subgroup of $G$).
Thus the vertex stabilizers of the action of $G$ on $\Sigma$
are all finite, and consequently this action is geometric.

\bigskip\noindent
{\bf 6.3 Systolic boundary of a systolic simplicial complex.}

\medskip
We recall, mostly from [OP], the necessary informations concerning
the concept of systolic boundary. 
For more informations, the reader is referred to
[OP] and to Subsection 9.3 in [OS].
Given a systolic simplicial complex 
$\Upsilon$, its systolic boundary $\partial\Upsilon$ is defined
using the objects called {\it good geodesic rays} 
(introduced in Definition 3.2 in [OP]).
For our purposes, we only need some properties of good
geodesic rays, which we recall below (see Lemmas 6.3.1 and 6.3.3), 
and here we only mention that
they are some special geodesic rays in the 1-skeleton of a systolic complex. 
As a set, {\it systolic boundary} $\partial\Upsilon$ is then
the set of all good geodesic rays in $\Upsilon$
quotiened by the equivalence relation of being at finite Hausdorff
distance from one another in $\Upsilon$ ([OP], Definition 3.6).

The following useful property follows immediately from Corollary 3.10
in [OP].

\medskip\noindent
{\bf 6.3.1 Lemma.}
{\it For any vertex $O$ in a systolic simplicial complex $\Upsilon$,
and any point $\xi\in\partial\Upsilon$, there is a good geodesic ray
$r$ in $\Upsilon$ started at $O$ and representing $\xi$.}

\medskip
The next result follows fairly directly from the definition of a good
geodesic ray (as given in [OP]) and from the structure of a tree
of systolic complexes. We do not present the details of a strightforward
proof of this result, but only an outline.

\medskip\noindent
{\bf 6.3.2 Lemma.} 
{\it Let $p:\Upsilon\to X$ be a tree of systolic complexes, and let $t$ be any
vertex in the tree $X$. Then any good geodesic ray in the systolic
complex $p^{-1}(t)$ is also a good geodesic ray in $\Upsilon$.} 

\medskip\noindent
{\bf Sketch of proof:}
We outline the straightforward argument which justifies the lemma,
referring the reader to [OP]  for explanations of the notions
appearing in this argument (which are used in that paper to define
good geodesic rays).

Obviously, the subcomplex $\Upsilon_t:=p^{-1}(t)$ is geodesically
convex in $\Upsilon$ (for the natural geodesic metric in the 1-skeleton).
Consequently, any directed geodesic in $\Upsilon_t$ is also a directed geodesic in $\Upsilon$. Furthermore, a surface spanned on
a loop in $\Upsilon_t^{(1)}$ is minimal in $\Upsilon_t$ iff it is minimal
in $\Upsilon$ (a surface in $\Upsilon$ spanned on such a loop
and not contained in $\Upsilon_t$ can be easily shown to be not minimal). 
It follows that any Euclidean geodesic in $\Upsilon_t$
is also a Euclidean geodesic in $\Upsilon$. In view of the definition
of a good geodesic ray (Definition 3.2 in [OP]), this completes the proof.

\medskip
Part (1) of the naxt lemma is a special case of Corollary 3.4 in [OP],
and part (2) coincides with Lemma 3.8 in the same paper.
Given a good geodesic ray $r$, we denote by $r(i):i\ge0$ its
consecutive vertices. We also denote by $d_{\Upsilon^{(1)}}$
the natural polygonal metric in the 1-skeleton of $\Upsilon$.

\medskip\noindent
{\bf 6.3.3 Lemma.}
{\it There is some universal constant $D>0$ 
satisfying the following properties.
For
any systolic complex $\Upsilon$ and any two good geodesic rays $r_1,r_2$ in $\Upsilon^{(1)}$
based at the same vertex $O$ we have:}
\item{(1)} {\it $d_{\Upsilon^{(1)}}(r_1(i),r_2(i))\le{i\over j}\cdot
d_{\Upsilon^{(1)}}(r_1(j),r_2(j))+D$
for any integer $i,j$ such that $0<i<j$;}
\item{(2)} 
{\it $r_1,r_2$ represent the same point in $\partial\Upsilon$ (i.e.\ they lie at
finite Hausdorff distance from one another) iff $d_{\Upsilon^{(1)}}(r_1(i),r_2(i))\leq D$ for all positive
integers $i$.}

\medskip
We now pass to describing the topology of $\partial\Upsilon$.
To do this, we fix a vertex $O$ in $\Upsilon$, and we denote by
${\cal R}_{O,\Upsilon}$ the set of all good geodesic rays in $\Upsilon$
started at $O$. Note that, in view of Lemma 6.3.1, this set contains
good geodesic rays representing all points of $\partial\Upsilon$.

Following Section 4 in [OP], 
the topology of $\partial\Upsilon$ is introduced by means
of local neighbourhood systems (which consist of sets that are
not necessarily open in the resulting topology).
More precisely, for each $\xi\in\partial\Upsilon$ we have a family ${\cal N}_\xi$
of sets containing $\xi$, called {\it standard neighborhoods} of $\xi$,
and the whole system ${\cal N}_\xi:\xi\in\partial\Upsilon$
satisfies some appropriate axioms. Open sets are described as those
$U\subset\partial\Upsilon$ for which $\forall \xi\in U$ $\exists Q\in{\cal N}_\xi$
such that $Q\subset U$. Moreover, each $Q\in{\cal N}_\xi$ contains some
open neighborhood of the point $\xi$.
Finally, standard neighborhoods $Q\in{\cal N}_\xi$ have the form
$$
Q=Q(r,N,R)=\{ \eta\in\partial\Upsilon:\hbox{ for some }
r'\in[\eta]\cap{\cal R}_{O,\Upsilon}\hbox{ it holds }
d_{\Upsilon^{(1)}}(r(N),r'(N))\leq R  \},
$$
for any good geodesic ray $r\in[\xi]\cap{\cal R}_{O,\Upsilon}$, 
and any positive integers $N,R$ with $R\geq D+1$,
where $D$ is a constant as in Lemma 6.3.3, and where $[\eta]$ denotes
the equivalence class of good geodesic rays representing the point $\eta$
(see Definition 4.1 in [OP]).

It is shown in [OP] that the above system of standard neighbourhoods satisfies
the appropriate axioms (Proposition 4.4), and that the resulting topology
in $\partial\Upsilon$ does not depend on the choice of a vertex $O$
(Lemma 5.5). 

Next result records basic properties of systolic boundaries
(see Corollaries 5.2 and 5.4, and Propositions 5.3 and 5.6 in [OP]).

\vfill\break

\medskip\noindent
{\bf 6.3.4 Lemma.}
\item{(1)}
{\it If a systolic complex $\Upsilon$ is uniformly
locally finite then its systolic boundary $\partial\Upsilon$ is compact, metrisable
and has finite topological dimension.}
\item{(2)}
{\it A locally finite systolic complex 
has non-empty systolic boundary  
iff it is unbounded.}

\medskip
Note that Lemma 6.3.4 applies in particular to systolic complexes acted upon
geometrically by a group.

A useful addition to the above description of topology in $\partial\Upsilon$ 
is the following
characterization of convergence.

\medskip\noindent
{\bf 6.3.5 Lemma.}
{\it Let $\xi$ and $\xi_n:n\ge1$ be points of $\partial\Upsilon$,
and let $r$ and $r_n:n\ge1$ be good geodesic rays in $\Upsilon$,
started at a fixed vertex $O$, representing these points, respectively.
Then the sequence $(\xi_n)$ is convergent to $\xi$ 
iff for some $R\ge D+1$ the sequence
$$
(\max\{ i:d_{\Upsilon^{(1)}}(r_n(i),r(i))\le R \})_{n\ge1}
$$
diverges to $+\infty$.}

\medskip\noindent
{\bf Proof:}
It follows from Lemma 4.3 in [OP] that for any 
good geodesic ray $r'$ in $\Upsilon$ started at $O$ and 
representing the same
point $\xi$ as $r$, to each $N',R'$ one can
associate $N$ such that $Q(r,N,R)\subset Q(r',N',R')$.
Thus, as a basis of neighbourhoods of
$\xi$ it is sufficient to take the family
$Q(r,N,R):N\ge1$. The lemma follows then directly from
the definition of standard neighbourhoods $Q$ given above.

\bigskip\noindent
{\bf 6.4 Proof of Theorem 0.3(4).}

\medskip
Let $\Sigma=\Sigma({\cal G},\Xi)$ be the systolic complex described 
in Subsection 6.2. We need to show that the systolic boundary
$\partial\Sigma$ is homeomorphic to the dense amalgam 
$\widetilde\sqcup(\partial\Sigma_1,\dots,\partial\Sigma_k)$.
Note that, under assumptions of the theorem, we obviously have 
that a group $G_{v_i}$ is finite iff the associated complex $\Sigma_i$
is bounded, and, by Lemma 6.3.3(2), this happens iff 
$\partial\Sigma_i=\emptyset$.

Similarly as in the proof of Lemma 4.4.1,
consider first the case when all groups $G_{v_i}$ are finite.
By our convention, we have that the dense amalgam
$\widetilde\sqcup(\partial\Sigma_1,\dots,\partial\Sigma_k)$
is then homeomorphic to the Cantor space $C$. 
By Lemma 4.1.6(2), 
the Bass-Serre tree $X$ is then infinite, locally finite and such
that the vertices which split $X$ into at least 3 unbounded components
form a net in $X$.
Moreover, the subcomplexes
$\Sigma_t=p^{-1}(t)$ for the natural structure of a tree of systolic complexes 
$p:\Sigma\to X$ are uniformly bounded
(because each such subcomplex is somorphic to one of the complexes
$\Sigma_i$).
It is not hard to observe that in such situation $p$ establishes the natural
bijective correspondence between the classes of good geodesic rays in 
$\Sigma$ and the ends of $X$, and that the systolic boundary 
$\partial\Sigma$ has then the natural topology of $\partial X$,
and hence it is homeomorphic to the Cantor space $C$.
We omit further details and conclude that the theorem follows
in the considered case.

In the remaining case, after possibly permuting the indices,
we have that for some $m\in\{ 1,\dots,k\}$
the vertex groups $G_{v_1},\dots,G_{v_m}$ are infinite,
while the remaining ones are finite.
Since the boundaries $\partial\Sigma_j$ for $j>m$ are then empty,
by our convention we have
$\widetilde\sqcup(\partial\Sigma_1,\dots,\partial\Sigma_m,\dots,
\partial\Sigma_k)=
\widetilde\sqcup(\partial\Sigma_1,\dots,\partial\Sigma_m)$.
Thus we need to show that 
$\partial\Sigma$ is homeomorphic to 
$\widetilde\sqcup(\partial\Sigma_1,\dots,\partial\Sigma_m)$.
In view of Theorem 0.2, and using the terminology introduced
at the beginning of Section 2, we need to show that $\partial\Sigma$ is
$(\partial\Sigma_1,\dots,\partial\Sigma_m)$-regular.

We start with describing a family 
${\cal Y}={\cal Y}_1\sqcup\dots\sqcup{\cal Y}_m$
of subsets in $\partial\Sigma$, and then we show that it satisfies
(the appropriate version of) conditions (a1)--(a5) of the introduction.
Recall that for each vertex $t$ of $X$ the subcomplex 
$\Sigma_t=p^{-1}(t)$ is systolic, and consider the map 
$\iota_t:\partial\Sigma_t\to\partial\Sigma$ defined as follows.
If $r$ is a good geodesic ray in $\Sigma_t$ representing a point
$\xi\in\partial\Sigma_t$, it follows from Lemma 6.3.2 that $r$ is also
a good geodesic ray in $\Sigma$. Thus, 
it represents a point $\eta\in\Sigma$, and we put $\iota_t(\xi):=\eta$.
Since $\Sigma_t$ is geodesically convex in $\Sigma$, the map $\iota_t$
is well defined and injective. By referring to the characterization
of convergence provided in Lemma 6.3.5, this also easily implies that
$\iota_t$ is continuous. Since by Lemma 6.3.4(1) the boundary $\partial\Sigma_t$ is compact, it follows that $\iota_t$ is an embedding.

Recall that $\pi:X\to Y$ is the canonical projection from the Bass-Serre 
tree to the underlying graph of the graph of groups $\cal G$.
For each $i\in \{ 1,2,\dots,m \}$ put
${\cal Y}_i:=\{ \iota_t(\partial\Sigma_t):t\in V_X, \pi(t)=v_i \}$.
We turn to verifying conditions (a1)--(a5).

\smallskip
It follows from the construction
of $\Sigma$ that
for each $t$ with $\pi(t)=v_i$,  the subcomplex $\Sigma_t$ is isomorphic to
$\Sigma_i$. Since each $\iota_t$ is an embedding, it follows that
the sets in each ${\cal Y}_i$ are all homeomorphic to 
$\partial\Sigma_i$. To complete verification
of (a1), we need to show that for distinct vertices $t,s\in V_X$
the images $\iota_t(\partial\Sigma_t),\iota_s(\partial\Sigma_s)$ 
are disjoint. Let $\xi_t,\xi_s$ be any points in the boundaries
$\partial\Sigma_t$ and $\partial\Sigma_s$, respectively.
Let $r_t,r_s$ be good geodesic rays in the complexes 
$\Sigma_t,\Sigma_s$ representing the points $\xi_t$ and $\xi_s$,
respectively. It is clear from the structure of a tree of systolic complexes
for $\Sigma$ provided by the projection $p$ that the Hausforff distance
between $r_t$ and $r_s$ in $\Sigma$ is infinite, and hence
the points $\iota_t(\xi_t),\iota_s(\xi_s)\in\partial\Sigma$ do not coincide.
This finishes the verification of (a1).

\smallskip
To check that the family $\cal Y$ is null (i.e. to verify condition (a2)), 
due to compactness of
$\partial\Sigma$, it is sufficient to show the
following property:
{\it let $\xi_i$ be a convergent sequence of points in $\partial\Sigma$
such that for each $i$ there is $t_i\in V_X$ with 
$\xi_i\in\iota_{t_i}(\partial\Sigma_{t_i})$, and such that the vertices
$t_i$ are pairwise distinct; then any other sequence of points
$\xi_i'$ such that $\xi_i'\in\iota(\Sigma_{t_i})$ is also convergent
in $\partial\Sigma$ (in fact, to the same limit point as 
the sequence $\xi_i$).} To verify the above property, 
suppose that $\xi_i\to\xi$.
Let $r$, $r_i$ and $r_i'$ be good geodesic rays
in $\Sigma$ started at a fixed vertex $O$ and representing
$\xi$, $\xi_i$ and $\xi_i'$, respectively.
Consider first the case when $\xi\in\iota_t(\partial\Sigma_t)$ 
for some $t\in V_X$.
It follows that, except for some finite initial part, $r$ is contained
in $\Sigma_t$. 
By the convergence criterion of Lemma 6.3.5 applied to 
the convergence $\xi_i\to\xi$,
$r_i$ intersects $\Sigma_t$ for all $i$ large enough.
Moreover, since $t_i$ are pairwise distinct, for all $i$ large enough
the rays $r_i$ exit the subcomplex $\Sigma_t$ after intersecting it. 
For such $i$, let $r_i(j_i)$ be the last vertex on $r_i$ belonging
to $\Sigma_t$, i.e. the vertex through which $r_i$ exits $\Sigma_t$.
By applying again Lemma 6.3.5 to the convergence $\xi_i\to\xi$,
we conclude that $j_i\to\infty$. Further,
since both $\xi_i$ and $\xi_i'$ belong to the same subset $\iota_{t_i}(\partial\Sigma_{t_i})$, it follows that for each $i$ large enough
the ray $r_i'$ also intersects $\Sigma_t$,
and then exits it through some vertex $x_i$ which is at distance
at most 1 from $r_i(j_i)$ (because both $r_i$ and $r_i'$ exit $\Sigma_t$
through the same simplex). By Lemma 6.3.3(1), we get that
$d_{\Sigma^{(1)}}(r(j),r'(j))\le D+1$ for all $j\le j_i$.
Consequently, applying again Lemma 6.3.5 and the fact that $j_i\to\infty$,
we conclude that
$\xi_i'\to\xi$. Thus the property above follows in this case.

In the remaining case we consider $\xi$ which does not belong
to any subset $\iota_t(\partial\Sigma_t)$. Consequently, the ray $r$
exits every subcomplex $\Sigma_t$ that it intersects.
Let $(t_n)_{n\ge1}$ be the vertices such that $\Sigma_{t_n}$
are the consecutive subcomplexes intersected by $r$.
For each $i$, let $n(i)$ be the largest $n$ such that the ray $r_i$
intersects $\Sigma_{t_n}$. Since $\xi_i\to\xi$, we deduce using Lemma
6.3.5 that $n(i)\to\infty$ as $i\to\infty$. By the fact that both $\xi_i$
and $\xi_i'$ belong to $\iota_{t_i}(\partial\Sigma_{t_i})$,
we get that for each $i$ the ray $r_i'$ also intersects $\Sigma_{t_{n(i)}}$.
Since $n(i)\to\infty$, it follows by applying once again Lemma 6.3.5
that $\xi_i'\to\xi$, which completes the verification of (a2).

\smallskip
To check condition (a3), we need to show that any subset
$Z=\iota_t(\partial\Sigma_t)$ is boundary in $\partial\Sigma$.
Fix any point $\xi\in\iota_t(\partial\Sigma_t)$, and let $r$ be a good
geodesic ray started at $O$ and representing $\xi$. As we have 
already noticed before, the ray $r$ is then contained in $\Sigma_t$,
except possibly some finite intial part.
We will construct a sequence $\xi_n$ of points in
$\partial\Sigma\setminus\iota_t(\partial\Sigma_t)$ which converges to $\xi$.
In the argument, we use the notation of Subsection 6.2.

Denoting $t=(gG_{v_i},v_i)$ for some $i\le m$, 
we get that the subcomplex $\Sigma_t$ is preserved by the subgroup
$gG_{v_i}g^{-1}<G=\pi_1({\cal G},\Xi)$. The restricted action of this
subgroup on $\Sigma_t$ is geometric (in particular, cocompact),
because it is equivariantly isomorphic to the action of $G_{v_i}$
on $\Sigma_i$. 
Let $a\in O_Y$ be an unoriented edge in the underlying graph $Y$
of $\cal G$, started at $v_i$, and satisfying the assertion of Lemma 4.1.6(1).
Denote by $E_{a,t}$ the set of all nonoriented edges of $X$ which are
the lifts of $|a|$ under the projection $\pi:X\to Y$, and which are adjacent
to $t$. Denote also by $S_{a,t}$ the set of all simplices in $\Sigma_t$
of form $\kappa_\varepsilon\cap\Sigma_t:\varepsilon\in E_{a,t}$.
The subgroup $gG_{v_i}g^{-1}$ obviously acts transitively on $S_{a,t}$,
and hence, by cocompactness of its action on $\Sigma_t$,
there is a constant $D_0$ such that
for each $n$ sufficiently large (namely for those $n$ for which
$r(n)\in\Sigma_t$) the vertex $r(n)$ lies at the distance
at most $D_0$ from some simplex in $S_{a,t}$, say
$\kappa_{\varepsilon_n}\cap\Sigma_t$. 
For each edge $\varepsilon_n$ as above, 
by referring to the assertion of Lemma 4.1.6(1),
choose a vertex $t_n$ in $X$ such that 
\item{$\bullet$} $\pi(t_n)=v_j$ for some $j\le m$ (so that $\Sigma_{t_n}$
is unbounded, as being isomorphic to $\Sigma_{j}$);
\item{$\bullet$}
the path from $t$ to $t_n$ in $X$ passes through the edge $\varepsilon_n$.

\noindent
Choose any point $\xi_n\in\iota_{t_n}(\partial\Sigma_{t_n})$
(which exists due to Lemma 6.3.4(2)),
and let $r_n$ be a good geodesic ray in $\Sigma$ started at $O$
and representing $\xi_n$. Note that, due to the structure of $\Sigma$
as the tree of systolic complexes, for $n$ sufficiently large the ray $r_n$
intersects $\Sigma_t$ and exits it through the simplex 
$\kappa_{\varepsilon_n}$. It follows that for those $n$ we have 
$\xi_n\notin\iota_t(\partial\Sigma_t)$, and that
$d_{\Sigma^{(1)}}(r(n),r_n(n))\le D_0+1$.  Applying Lemma 6.3.5, 
we deduce from the latter that $\xi_n\to\xi$,
hence condition (a3).

\smallskip
To check (a4), suppose that 
$\xi\in\partial\Sigma\setminus\bigcup_{t\in V_X}\iota_t(\partial\Sigma_t)$.
We need to find a sequence of points 
$\xi_n\in\bigcup_{t\in V_X}\iota_t(\partial\Sigma_t)$ converging to $\xi$.
Let $r$ be a good geodesic ray in $\Sigma$ started at $O$ and
representing $\xi$. By the above assumption on $\xi$, the ray
$r$ exits every subcomplex $\Sigma_t$ that it intersects.
Let $(t_n)_{n\ge0}$ be the sequence of vertices in $X$ such that
$\Sigma_{t_n}$ are the consecutive subcomplexes of this form 
intersected by $r$.
Choose any integers $j_n$ such that $r(j_n)\in\Sigma_{t_n}$,
and note that $j_n\to\infty$. Denote also by $\kappa_{\varepsilon_n}$
the simplex through which the ray $r$ eaxits $\Sigma_{t_n}$.
Since $r$ projects through $p$ on an infinite ray $\rho$ in $X$ 
(formed of the 
consecutive edges $\varepsilon_n$), for each $n$ there is a vertex
$s_n$ in $X$, with $\pi(s_n)=v_i$ for some $i\le m$, and such that
the path from $t_0$ to $s_n$ shares first $n+1$ edges with $\rho$.
Let $\xi_n$ be any point in $\iota_{s_n}(\partial\Sigma_{s_n})$
(which exists by Lemma 6.3.4(2)), and let $r_n$ be a good geodesic
ray in $\Sigma$ started at $O$ and representing $\xi_n$. By the choice
of $s_n$ and $\xi_n$, the ray $r_n$ intersects $\Sigma_{t_n}$
and exits it through the simplex $\kappa_{\varepsilon_n}$.
It follows that for some $j\ge j_n$ we have 
$d_{\Sigma^{(1)}}(r(j),r_n(j))\le1$. Since $j_n\to\infty$,
we deduce from Lemma 6.3.5 that $\xi_n\to\xi$,
which completes the verification of (a4).

\smallskip
It remains to check condition (a5).
Let $\xi_1,\xi_2\in\partial\Sigma$ be any two points which do not
belong to the same subset of $\cal Y$. Let $r_1,r_2$ be any good
geodesic rays in $\Sigma$ started at $O$ and representing
$\xi_1$ and $\xi_2$, respectively. By the above assumption on
$\xi_1,\xi_2$, up to transposing the indices,
there is an edge $\varepsilon$ in $X$ such that
the path $p\circ r_1$ in $X$, being the projection of $r_1$ to $X$,
passes through $\varepsilon$, while the path $p\circ r_2$ does not.
We define a splitting of $\partial\Sigma$ into two subsets.
Let $H_1$ consist of all points of $\partial\Sigma$ which can be represented
by a good geodesic ray $r$ started at $O$ and such that $p\circ r$
passes through $\varepsilon$. 
Put also $H_2=\partial\Sigma\setminus H_1$, and note that it
consists of all points of $\partial\Sigma$ which can be represented
by a good geodesic ray $r$ started at $O$ and such that $p\circ r$
does not pass through $\varepsilon$. It is an easy observation
involving Lemma 6.3.5 that both $H_1$, $H_2$ are closed in
$\partial\Sigma$. Obviously, they are also $\cal Y$-saturated,
and separate $\xi_1$ from $\xi_2$. This justifies condition (a5),
and completes the proof of Theorem 0.3(4).

\bigskip
\noindent
{\bf 7. 
Boundaries of Coxeter groups.}

\medskip
In this section we  prove Theorem 0.4 of the introduction. 
The rough idea of the proof is this.
First, we observe (by referring to the characterization of dense amalgams)
that any non-elementary splitting of $W$ as free product of special subgroups,
amalgamated along a finite special subgroup, leads to the expression of
the boundary of $W$ as the dense amalgam of boundaries of the factors
(see Proposition 7.3.1 for precise statement).
Next, we note that splittings as above correspond to splittings of the nerve
of $W$ along separating simplices (including the empty one).
Further, inspired by the comments in Section 8.8 of [Da],
we argue that on the level of groups the terminal factors
of iterations of such splittings are either the maximal finite or the maximal 
1-ended special subgroups
(this is a more precise version of the assertion of Proposition 8.8.2 in [Da]).
Finally, applying general properties 
of the operation of dense amalgam (given in Proposition 0.1),
we show that this yields the assertion. 
Details are provided in Subsections 7.1--7.4 below.

\bigskip\noindent
{\bf 7.1 Decompositions of simplicial complexes.}

\smallskip
We introduce a usefull terminology, and provide basic facts, concerning 
decompositions 
of simplicial complexes along simplices. 
The idea of such decompositions is well known in graph theory,
see e.g. [Di].
Since we need only very basic facts, and in a specific setting,
we briefly provide an independent account.

\medskip\noindent
{\bf 7.1.1 Definition.}
Let $L$ be a simplicial complex. A {\it splitting} of $L$ along a simplex
is an expression of $L$ as the union of proper nonepty subcomplexes
$L_1,L_2$ whose intersection $L_1\cap L_2$ is either empty
or a single simplex. $L_1$ and $L_2$ are then called 
{\it the parts} of this
splitting. A simplicial complex $L$ is {\it irreducible} if it has no splitting.

\medskip
Observe that the parts of any splitting of $L$ are 
full subcomplexes of $L$. Note also 
that $L$ is irreducible if it is connected and has no separating simplex.

We now give a recursive definition of a decomposition of a simplicial
complex, and of its factors.

\medskip\noindent
{\bf 7.1.2 Definition.}
 {\it A decomposition} of a simplicial complex $L$ is any sequence
of splittings from the following recursively described family:
\item{$\bullet$} the empty sequence of splittings forms the trivial
decomposition of $L$, and the set of factors of this decomposition
is $\{ L \}$;
\item{$\bullet$} a single splitting of $L$ along a simplex is a decomposition,
and its set of factors is the set of two parts of the splitting;

\item{$\bullet$}
if a sequence of splittings is a decomposition of $L$, and if 
$\{ L_1,\dots,L_m \}$ is the set of its factors, then adding to this sequence
a splitting of one of those factors, say $L_m$, we also get a decomposition of $L$;
moreover, if $L_m',L_m''$ are the parts of the above splitting of $L_m$,
the set of factors of the new decomposition is $\{ L_1,\dots,L_{m-1},L_m',L_m'' \}$.

\medskip
Note that it may happen that $L_m'$ or $L_m''$ as above coincides with $L_j$
for some $j\le m-1$, but then of course this subcomplex appears just once in the
set of factors of the corresponding decomposition.

\smallskip
A decomposition of $L$ is {\it terminal} if its every factor is irreducible. Obviously,
every finite simplicial complex admits a terminal decomposition. Next lemma shows
that any two terminal decompositions of a finite simplicial complex share the sets of
factors (though they may be quite different as sequences of splittings).
Thus, we call the factors of any terminal decomposition as above
the {\it terminal factors}.
The same lemma characterizes the terminal factors of a finite
simplicial complex.
In its statement we use the term {\it maximally full irreducible} subcomplex, which
denotes any subcomplex which is maximal for the inclusion in the family of all
full and irreducible subcomplexes of a given complex.

\medskip\noindent
{\bf 7.1.3 Lemma.}
{\it The set of factors of any terminal decomposition of a finite simplicial complex $L$
coincides with the set of all maximally full irreducible subcomplexes of $L$.}
 
\medskip\noindent
{\bf Proof:} We start with showing two auxilliary claims.

\smallskip\noindent
{\bf Claim 1.} 
{\it Any factor of a terminal decompostion of $L$ is a maximally full irreducible
subcomplex of $L$.}

\smallskip
To prove Claim 1, consider any factor $K$ of a fixed terminal decomposition of $L$.
$K$ is clearly full and irreducible. Suppose, by contradiction, that $K$ is not maximally
full irreducible, and let $M$ be a full and irreducible subcomplex of $L$ containing $K$
as proper subcomplex. Denote by $\cal F$ the set of factors of our fixed terminal decomposition
of $L$. This decomposition induces a decomposition of $M$ for which the set of factors
is $\{ A\cap M:A\in{\cal F}, A\cap M\ne\emptyset \}$. 
In particular, $K\cap M=K$
is a factor of this induced decomposition of $M$, which contradicts irreducibility of $M$,
thus completing the proof of Claim 1.

\smallskip\noindent
{\bf Claim 2.}
{\it Given any decomposition of $L$, every irreducible subcomplex of $L$ is contained in
at least one factor of this decomposition.}

\smallskip
To prove Claim 2, note that if $M$ is an irreducible subcomplex of $L$, and if $L_1,L_2$
are the parts of some splitting of $L$ along a simplex, then $M\subset L_1$ or
$M\subset L_2$. Claim 2 then follows by iterating this observation.

\smallskip
Now, in view of Claim 1, to prove Lemma 7.1.3, it is sufficient to show
that any maximally full irreducible subcomplex $M$ of $L$ is a factor
in every terminal decomposition of $L$. Fixing such a decomposition,
we get from Claim 2 that $M$ is contained in some factor $K$ of this
decomposition. Since $K$ is full and irreducible, maximality of $M$
implies that $M=K$, which completes the proof.

\medskip\noindent
{\bf 7.1.4 Example.}
The class of finite simplicial complexes in which all terminal factors
are simplices is well known. It coincides with the class of finite flag
simplicial complexes which contain no full subcomplex isomorphic
to a triangulation of the circle $S^1$, see [Dir]. According to the
terminology from [JS], which we follow, such complexes are called
$\infty$-{\it large}. In Section 8.8 of [Da], such complexes are called 
(a bit informally) {\it trees of simplices}. 1-skeletons of such 
complexes are known in graph theory as {\it chordal graphs}.

\bigskip\noindent
{\bf 7.2 Nerves of Coxeter systems.}

\medskip
Recall that the {\it nerve} $L=L(W,S)$ of a Coxeter
system $(W,S)$ is the simplicial complex whose vertex set coincides
with $S$, and whose simplices correspond to those subsets $T\subset S$
which span finite special subgroups $W_T<W$. In this subsection
we recall from [Da] few results and observations concerning properties of groups 
$W$ that can be read from properties of their nerves. 
The first fact below is straightforward
(compare [Da], Proposition 8.8.1).

\medskip\noindent
{\bf 7.2.1 Lemma.} 
{\it Let $(W,S)$ be a Coxeter system
with the nerve $L$, and let $S_1,S_2$ be the vertex sets of the parts 
$L_1,L_2$ of some splitting of $L$ along a simplex. For $i=1,2$
denote by $W_i$ the special subgroup of $W$ generated by $S_i$,
and denote by $W_0$ the special subgroup generated by 
the intersection $S_1\cap S_2$
(in particular, the trivial subgroup if $S_1\cap S_2=\emptyset$). Then $W=W_1*_{W_0}W_2$, i.e. $W$
is the free product of the subgroups $W_1,W_2$ amalgamated
along the finite subgroup $W_0$.}

\medskip\noindent
{\bf 7.2.2 Theorem} ([Da], Theorem 8.7.2).
{\it A Coxeter group $W$ is 1-ended iff its nerve is an irreducible
simplicial complex distinct from a simplex.}

\medskip
Note that the groups appearing in Theorem 7.2.2 are precisely
those Coxeter groups which have nonempty connected boundary.
(This follows e.g. from Proposition 8.6.2(i) and Theoerm I.8.3(ii)
in [Da].)

\medskip\noindent
{\bf 7.2.3 Theorem} ([Da], Theorem 8.7.3).
{\it A Coxeter group $W$ is 2-ended iff it can be expressed as 
the product $W=W_0\times W_1$, where $W_0$ is a special subgroups
isomorphic to the infinite dihedral group, and $W_1$ is a finite special subgroup (including the case of the trivial subgroup).}

\medskip
Note that the groups appearing in Theorem 7.2.3 are precisely those
Coxeter groups whose boundaries are the spaces consisting 
of two points. Moreover, nerves of such groups are
suspended simplices (including the case of the suspended empty simplex).
However, not every Coxeter group whose nerve is a suspended simplex
is 2-ended.

\medskip\noindent
{\bf 7.2.4 Proposition} ([Da], Proposition 8.8.5).
{\it A Coxeter group is virtually free nonabelian
iff it is not 2-ended and its nerve is an $\infty$-large simplicial complex
distinct from a simplex.}

\medskip
Note that the groups appearing in Proposition 7.2.4 all have
Cantor space $C$ as the boundary. In fact, it is not difficult to show
(and it follows in particular from Theorem 0.4)
that the condition in the proposition fully characterizes
the Coxeter groups which have Cantor space $C$ as their boundaries.

\vfill\break

\bigskip\noindent
{\bf 7.3 Dense amalgams and decompositions of nerves.}

\medskip
We start with the basic observation, Proposition 7.3.1 below,
bringing dense amalgams into considerations concerning
boundaries of Coxeter groups. Since the proof of this proposition
goes along the same lines as the proof of Theorem 0.3(4)
given in Subsection 6.4, we omit it. We only note that, in view of Lemma 7.2.1,
a splitting of the nerve of $W$ along a simplex induces a splitting
of $W$ over a finite group; moreover, the assumption below concerning indices means that the corresponding graph of groups of the splitting
is non-elementary. This makes Proposition 7.3.1 completely analogous
to the results in parts (1)--(4) of Theorem 0.3 (or rather to their special cases,
with $\cal G$ corresponding to a single amalgamated free product). 

\medskip\noindent
{\bf 7.3.1 Proposition.}
{\it Under assumptions and notation of Lemma 7.2.1,
suppose additionally that for at least one of the indices $i\in\{ 1,2 \}$
we have $[W_i:W_0]\ne2$ (i.e. the subgroup $W_0$ has index
greater than 2 in at least one of the groups $W_i$).
Then 
$$\partial(W,S)\cong
\widetilde\sqcup\big( \partial(W_1,S_1),\partial(W_2,S_2) \big).$$}

\noindent
{\bf Remark.}
Note that if in the setting of Proposition 7.3.1 we have 
$[W_i:W_0]=2$ for both $i=1,2$ then 
$W_1\cong W_0\times Z_2\cong W_2$ and
$W\cong W_0\times D_\infty$. Then we obviously have
$\partial(W_1,S_1)\cong\partial(W_2,S_2)\cong\partial(W_0,S_1\cap S_2)$,
while the boundary $\partial W$ is homeomorphic to the suspension
of those spaces. This shows that the assumption in the proposition 
concerning indices
$[W_i:W_0]$ is essential.

\medskip
Next result is an extension of Proposition 7.3.1 to more complicated decompositions
of the nerves of Coxeter systems.

\medskip\noindent
{\bf 7.3.2 Proposition.}
{\it Suppose that $L_1,\dots,L_k$ are the factors of a decomposition
of the nerve $L$ of a Coxeter system $(W,S)$, and let $(W_i,S_i)$
be the Coxeter systems of special subgroups of $W$ corresponding
to the vertex sets $S_i$ of the subcomplexes $L_i$. Suppose also
that $W$ is not 2-ended, and that $k\ge2$. Then
$$
\partial(W,S)\cong
\widetilde\sqcup\big( \partial(W_1,S_1),\dots,\partial(W_k,S_k) \big).
$$}
{\bf Proof:}
We argue by induction with respect to the length $n$ of a sequence
of splittings along simplices that constitutes a decomposition of $L$
under consideration. Since we assume that the number $k$ of factors
is at least 2, we have $n\ge1$. The case $n=1$ follows by
Proposition 7.3.1. Thus, it remains to verify the general inductive step.

Suppose that the statement holds true for some decomposition of
length $n$, and that $L_1,\dots,L_k$ are the factors of this decomposition.
Consider a decomposition of length $n+1$ obtained by adding a splitting
of the factor $L_k$, with parts $L_k',L_k''$. 
Obviously, the set of factors of the new decomposition is then
$\{ L_1,\dots,L_{k-1},L_k',L_k'' \}$.
Denote by $S_k',S_k''$
the vertex sets of $L_k'$ and $L_k''$, respectively, and let $W_k',W_k''$
be the special subgroups generated by these sets.
We need to consider two cases.

\medskip\noindent
{\it Case 1:} at least one of the indices $[W_k':(W_k'\cap W_k'')]$ and
$[W_k'':(W_k'\cap W_k'')]$ is distinct from 2.

In this case the splitting of $L_k$ into $L_k'$ and $L_k''$ satisfies
the assumptions of Proposition 7.3.1, and hence
$\partial(W_k,S_k)\cong\widetilde\sqcup\big( 
\partial(W_k',S_k'),\partial(W_k'',S_k'') \big)$. Consequently,
using the inductive assumption and Proposition 0.1(2), we get 
$$
\partial(W,S)\cong
\widetilde\sqcup\big( \partial(W_1,S_1),\dots\partial(W_k,S_k) \big)
\cong
$$
$$
\cong
\widetilde\sqcup
\big[
\partial(W_1,S_1),\dots\partial(W_{k-1},S_{k-1}),
\widetilde\sqcup\big( 
\partial(W_k',S_k'),\partial(W_k'',S_k'') \big)  
\big]
\cong
$$
$$
\cong
\widetilde\sqcup
\big(  
\partial(W_1,S_1),\dots\partial(W_{k-1},S_{k-1}),
\partial(W_k',S_k'),\partial(W_k'',S_k'') 
\big).
$$
Now, if $L_k'$ or $L_k''$ coincides with one of the subcomplexes
$L_1,\dots,L_{k-1}$, we apply Proposition 0.1(3) to get the assertion.
Otherwise, the assertion follows directly.

\medskip\noindent
{\it Case 2:} $[W_k':(W_k'\cap W_k'')]=[W_k'':(W_k'\cap W_k'')]=2$.

Note that, under this assumption, the group $W_k$
is 2-ended, while both $W_k'$ and $W_k''$ are finite. Consequently,
the boundary $\partial(W_k,S_k)$ is the space consisting of 2 elements,
which we denote $Q_2$. We also have
$\partial(W_k',S_k')=\partial(W_k'',S_k'')=\emptyset$.
Using this, the inductive assumption, Proposition 0.1(4), and the 
properties of dense amalgam involving the empty set,
we get
$$
\partial(W,S)\cong
\widetilde\sqcup\big( \partial(W_1,S_1),\dots\partial(W_k,S_k) \big)
\cong
\widetilde\sqcup
\big(  
\partial(W_1,S_1),\dots\partial(W_{k-1},S_{k-1}),Q_2
\big)
\cong
$$
$$
\cong
\widetilde\sqcup
\big(  
\partial(W_1,S_1),\dots\partial(W_{k-1},S_{k-1})
\big)
\cong
\widetilde\sqcup
\big(  
\partial(W_1,S_1),\dots\partial(W_{k-1},S_{k-1}),\emptyset
\big).
$$
This implies the assertion, no matter if some of the boundaries
$\partial(W_j,S_j):1\le j\le k-1$ is empty or not.

This completes the proof.


\bigskip\noindent
{\bf 7.4 Proof of Theorem 0.4.}

\medskip
First, observe that the nerve $L$ of $(W,S)$ is not an $\infty$-large
simplicial complex. Indeed, it is not a simplex since $W$ is not finite.
It is not any other $\infty$-large simplicial complex by Theorem 7.2.3
and Proposition 7.2.4. In view of a comment provided in Example 7.1.4,
it follows from Lemma 7.1.3 that $L$ contains at least one
maximally full irreducible subcomplex distinct from a simplex.
Applying Theorem 7.2.2, this means that $W$ contains at least one 
maximal 1-ended special subgroup. Hence, we have shown the assertion
that $k\ge1$.

Now, consider any terminal decomposition of the nerve $L$,
and let $L_1,\dots,L_m$ be the factors of this decomposition.
For $i=1,\dots,m$, denote by $S_i\subset S$ the vertex set of $L_i$,
and by $W_i$ the special sungroup generated by $S_i$.
By Lemma 7.3.2, we get that
$$
\partial(W,S)\cong\widetilde\sqcup\big( \partial(W_1,S_1),\dots,
\partial(W_m,S_m) \big).
$$ 
Without loss of generality, suppose that $L_1,\dots,L_k$ are precisely
those factors among $L_1,\dots,L_m$ which are not simplices.
Then $W_1,\dots,W_k$ is the family of all maximal 1-ended special subgroups of $W$.
We also know that $k\ge1$. 

Since for $k+1\le j\le m$ the subcomplexes $L_j$ are simplices, the corresponding special subgroups $W_j$ are finite, and their boundaries
$\partial(W_j,S_j)$ are empty. Since adding the empty set to the list
of densely amalgamated spaces does not affect the result,
it follows that
$$
\partial(W,S)\cong\widetilde\sqcup\big( \partial(W_1,S_1),\dots,
\partial(W_k,S_k) \big),
$$ 
which finishes the proof.

\bigskip
\centerline{\bf References}
\medskip

\itemitem{[Be]} M. Bestvina, {\it Local homology properties of boundaries of groups}, Michigan Math. J. 43 (1) (1996), 123--139.

\itemitem{[BH]} M. Bridson, A. Haefliger, 
Metric Spaces of Non-Positive
Curvature, Grundlehren der mathematischen Wissenschaften 319, Springer,
1999.

\itemitem{[CP]} G. Carlsson, E. Pedersen, 
{\it Controlled algebra and the {N}ovikov conjectures for {$K$}- 
and {$L$}-theory},
Topology 34 (1995), 731--758. 

\itemitem{[ChO]} V. Chepoi, D. Osajda,
{\it Dismantlability of weakly systolic complexes and applications},
 	arXiv:0910.5444, to appear in Trans. AMS.

\itemitem{[Dav]} R. Daverman, Decompositions of Manifolds, 
Academic Press, 1986.

\itemitem{[Da]} M. Davis, {\it The geometry and topology 
of Coxeter groups}, {London Mathematical Society Monographs Series},
{vol. 32}, {Princeton University Press}, {Princeton}, {2008}.

\itemitem{[Di]} R. Diestel, {\it Simplicial decompositions of graphs -
some uniqueness results}, Journal of Combinatorial Theory, Series B 42
(1987), 133--145. 

\itemitem{[Dir]} G.A. Dirac, {\it On rigid circuit graphs},
Abh. Math. Sem. Univ. Hamburg 38 (1961), 71--76.

\itemitem{[Dra]} A. Dranishnikov, 
{\it On {B}estvina-{M}ess formula},
in "Topological and asymptotic aspects of group theory",
Contemp. Math., vol. 394, p. 77--85, Amer. Math. Soc., Providence, 2006.

\itemitem{[FL]} F. T. Farrell, J.-F. Lafont, 
{\it E{Z}-structures and topological applications,}
Comment. Math. Helv. 80 (2005), 103--121.

\itemitem{[JS]} T. Januszkiewicz, J. \'Swi\c{a}tkowski,
{\it Simplicial nonpositive curvature},
Publ. Math. IHES 104 (1) (2006), 1--85.

\itemitem{[Ma]} A. Martin,
{\it Non-positively curved complexes of groups and boundaries}, Geometry \& Topology 18 (2014), 31--102


\itemitem{[MS]}
{A. Martin, J.  \'Swi\c{a}tkowski},
{\it Infinitely ended hyperbolic groups with homeomorphic Gromov boundaries},
J. Group Theory, to appear,
 arXiv:1303.6774.

\itemitem{[OP]} D. Osajda, P. Przytycki, {\it Boundaries od systolic groups,} 
Geometry \& Topology 13 (2009), 2807--2880.

\itemitem{[OS]} D. Osajda,  J.  \'Swi\c{a}tkowski, 
{\it On asymptotically hereditarily aspherical groups,}
 \itemitem{} arXiv:1304.7651.

\itemitem{[P]} P. Przytycki,
{\it The fixed point theorem for simplicial nonpositive curvature,} Mathematical Proceedings of Cambridge Philosophical Society 144 (2008), 683--695. 

\itemitem{[Ros]} D. Rosenthal, 
{\it Continuous control and the algebraic {$L$}-theory assembly map},
Forum Math. 18 (2006), 193--209.

\itemitem{[Se]} J-P. Serre, Trees, Springer, 1980.

\itemitem{[Ti]}
 {C. Tirel},
   {\it Z-structures on products},
  Algebr. Geom. Topol. 11 (2011), 2587--2625.

\bye